\documentclass{article}

\newtheorem{fed}{\textbf{Definition}}[section]
\newtheorem{thm}[fed]{\textbf{Theorem}}
\newtheorem{lemma}[fed]{\textbf{Lemma}}
\newtheorem{ex}[fed]{\textbf{Example}}
\newtheorem{rem}[fed]{\textbf{Remark}}
\newtheorem{prop}[fed]{\textbf{Proposition}}
\newtheorem{cor}[fed]{\textbf{Corollary}}

\usepackage{amssymb,amsmath}
\usepackage{graphicx}
\usepackage{dsfont}
\usepackage{color}
\usepackage{amscd,amssymb,
graphicx,epsfig,
latexsym,amsmath}
\usepackage{amsmath}
\usepackage{verbatim}
\newcommand{\N}{\mathbb{N}}

\newcommand{\Z}{\mathbb{Z}}
\newcommand{\R}{\mathbb{R}}
\newcommand{\C}{\mathbb{C}}

\newcommand{\la}{\langle}
\newcommand{\ra}{\rangle}
\newcommand{\p}{\partial}
\newcommand{\id}{\text{id}}
\newcommand{\wt}{\widetilde}
\newcommand{\om}{\omega}
\newcommand{\eps}{\varepsilon}
\newcommand{\J}{\mathcal{J}}

\numberwithin{equation}{section}

\begin{document}
\parindent=0pt
\parskip=4pt

\title{Periodic orbits in the restricted three-body problem and Arnold's $J^+$-invariant}
\author{Kai Cieliebak, Urs Frauenfelder and Otto van Koert}
\maketitle
\tableofcontents

%%%%%%%%%%%%%%%%%%%%%%%%%%%%%%%%%%%%%%%%%%%%%%%%%%%%%%%%%%%%
\section{Introduction}
%%%%%%%%%%%%%%%%%%%%%%%%%%%%%%%%%%%%%%%%%%%%%%%%%%%%%%%%%%%%

By the Whitney-Graustein theorem, the rotation number of an immersed loop in the plane completely classifies such loops up to regular homotopy. During a generic regular homotopy three disasters can occur, namely triple intersections as well as direct and inverse self-tangencies. Arnold \cite{arnold} noticed that the theory of immersions in the plane becomes extremely rich if one allows only generic homotopies avoiding the three disasters. He introduced three invariants. The one of interest in this paper is the invariant $J^+$, which
is invariant under triple intersections and inverse self-tangencies but is sensitive to direct self-tangencies. 

Immersions without direct self-tangencies are of geometric interest because their conormal lifts give Legendrian embeddings in the unit cotangent bundle of the plane \cite{arnold}. They are also natural from a physical point of view: orbits of a particle on the plane moving in a conservative force field cannot have direct self-tangencies, whereas inverse self-tangencies can occur if the force is allowed to be velocity dependent (like the Coriolis force or the Lorenz force). 

To analyze this phenomenon we consider conservative force fields on the plane generated by a vector potential and a scalar potential, which are smooth except for a $1/r$ singularity at the origin for the scalar potential. Such systems describe the motion of an electron in the presence of a proton and additional electric and magnetic fields, which is why we call them {\em Stark-Zeeman systems}. On the other hand, Newton's law of gravitation is mathematically equivalent to Coulomb's law and the Coriolis force is modelled in the same way as the Lorentz force. Therefore, Stark-Zeeman systems also describe many interesting cases of the restricted three-body problem in celestial mechanics, in which a satellite experiences the gravitational forces of two primary masses (which might be the earth and the moon) as well as the centrifugal force and the Coriolis force. 

In a generic homotopy of periodic orbits of Stark-Zeeman systems, apart from triple intersections and inverse self-tangencies the following two events can occur. First, the orbit might have velocity zero for some isolated moment in time. If a family of periodic orbits moves through this event an exterior loop is created or destroyed as shown in Figure~\ref{fig:cusp-zerovelocity}. This phenomenon was first observed by Hill in his discovery of the "moon of maximal lunarity" and will be recalled in Section~\ref{families}. Secondly, the orbit might pass through the origin, corresponding to a collision of the satellite with one of the primaries in the three-body problem. If this occurs in a family a loop around the origin is created or destroyed as shown in Figure~\ref{fig:cusp-collision}. 
See also Figure~\ref{fig:homotopy} for a family of periodic orbits in the restricted $3$-body problem in which both events occur.
We refer to such homotopies as \emph{Stark-Zeeman homotopies}. The task we address is the now following.\\ \\
\textbf{Task: } \emph{Find invariants under Stark-Zeeman homotopies.}
\\ \\
We introduce two such invariants both based on Arnold's $J^+$-invariant. The first invariant $\J_1$
is a combination of Arnold's $J^+$-invariant with the winding number around the origin. For the second invariant $\J_2$ we Levi-Civita regularize collisions and then take the $J^+$-invariant of the regularized orbit. Our first result is
\\ \\
\textbf{Theorem A: } \emph{$\mathcal{J}_1$ and $\mathcal{J}_2$ are invariant under Stark-Zeeman homotopies.}
\\ \\
We then ask how independent the two invariants are and prove
\\ \\
\textbf{Theorem B: } \emph{If the winding number of the loop around the origin is odd the two 
invariants determine each other, while in the case of even winding number they are completely independent. }
\\ \\
This paper is organized as follows. 
In Section~\ref{families} we recall the role of families of orbits in the search for periodic orbits in celestial mechanics. 
In Section~\ref{stark-zeeman} we introduce the class of Stark-Zeeman systems and describe prominent examples (which include several integrable systems) as well as their Levi-Civita and Moser regularizations.
In Section~\ref{topology} we introduce the notion of a Stark-Zeeman homotopy. After recalling the definition of Arnold's $J^+$-invariant, we introduce the invariants $\mathcal{J}_1$ and $\mathcal{J}_2$ and prove Theorems A and B. 

{\bf Acknowledgements. }
K.C. was supported by DFG grant CI 45/8-1, 
U.F. by DFG grant FR 2637/2-1, and
O.v.K. by NRF grant NRF-2016R1C1B2007662. 

%%%%%%%%%%%%%%%%%%%%%%%%%%%%%%%%%%%%%%%%%%%%%%%%%%%%%%%%%%%%
\section{Families of periodic orbits in the restricted three-body problem}\label{families}
%%%%%%%%%%%%%%%%%%%%%%%%%%%%%%%%%%%%%%%%%%%%%%%%%%%%%%%%%%%%

According to Poincar\'e's famous dictum, ``periodic orbits are the only tool to enter an otherwise impenetrable stronghold'' \cite{poincare}. Therefore the search for periodic orbits in the restricted $3$-body problem has a long history. 
Everybody having some experience with the search of these objects knows that one does not find them by trial and error. Instead, one usually looks at families of periodic orbits. The family parameter can be for example the energy or the mass ratio of the two primaries. One starts the family with an obvious periodic orbit. This can for example be a critical point, i.e., a Lagrange point, or a periodic orbit in a completely integrable system. For instance, as the energy goes to minus infinity the
restricted $3$-body problem just becomes the Kepler problem, for which on each negative energy level one has
a circular periodic orbit. As this orbit can be traversed forwards and backwards it gives rise to two
periodic orbits in the restricted $3$-body problem, the retrograde and the direct one. Strictly speaking, there are actually four of these orbits since one has direct and retrograde orbits around each of the primaries. Starting with these orbits for very low energy one then carries out a numerical homotopy by increasing the energy
to obtain a family of periodic orbits.

The pioneer of this approach was Hill \cite{hill}. In his theory of the moon he considered a limiting case
of the restricted $3$-body problem in which one of the primaries (the sun) becomes infinitely much heavier
than the second primary but in turn is infinitely far away. Hill's major interest was the direct periodic orbit,
since this orbit provides a good model for the orbit of the moon around the earth. In the actual restricted $3$-body problem, pioneers of this method were Darwin \cite{darwin1, darwin2} and Moulton \cite{moulton}. In particular, Darwin considered the case where one of the primaries is ten times heavier than the other. 
An extremely detailed study over many years was carried out by Str\"omgren and his school at the observatory of Copenhagen \cite{stromgren}. The case considered by Str\"omgren is where the two primaries have equal mass.
The current notation of these families is as well due to the school of Str\"omgren. In this notation the family of the retrograde orbit is referred to as the $f$-family, the family of the direct one as the $g$-family, where the families $a$ through $e$ start at the Lagrange points. For more information we refer to the book by Szebehely \cite{szebehely}. 

A different approach to Hill's direct periodic orbit was found by Ljapunov \cite{ljapunov}. In Hill's lunar
problem three forces are acting: the gravitational force of the earth, the Coriolis force, and a third force
which one can interpret as a tidal force. This tidal force arises in the limit process in which the gravitational force
of the infinitely heavy but infinitely far away sun and the centrifugal force centered in the infinitely far
away sun cancel each other to first order. While the gravitational force of the earth and the Coriolis force are
rotationally invariant, the tidal force breaks the rotational symmetry. Ljapunov symmetrizes the tidal force and
then considers a homotopy between the tidal force and its symmetrization. During the homotopy he fixes the period so that he obtains a one-parameter family.  For the symmetrized tidal force
there exist circular periodic orbits. The homotopy leads to a family of periodic orbits which interpolates
between a circular periodic orbit in the symmetrized problem and the direct periodic orbit in Hill's lunar problem, see also the paper by Rjabov \cite{rjabov}.
 
Poincar\'e \cite{poincare} found another method to obtain periodic orbits in the restricted $3$-body problem, namely by fixing energy but varying the mass ratio of the two primaries. As the mass ratio goes to zero one
obtains the rotating Kepler problem, i.e., the Kepler problem in rotating coordinates. This problem is again completely integrable and its periodic orbits are well known, so that they can be used as the starting points for
a family. There are two kinds of periodic orbits in the rotating Kepler problem: circular ones and honestly elliptical ones. Poincar\'e referred to periodic orbits in the restricted three-body problem obtained from the circular ones as periodic orbits of the first kind, and to the latter ones as periodic orbits of the second kind. Deep interest in the periodic orbits of the second kind was sparkled during the Apollo program through the work of Arenstorf
\cite{abraham-marsden, arenstorf}. For more recent developments about this approach we refer to the books by Bruno \cite{bruno} and H\'enon \cite{henon}. 

In a family of periodic orbits some disasters can occur. The first disaster is that the velocity of a periodic orbit can vanish at some moment. As the acceleration at this moment need not vanish, the periodic orbit need not be constant. However, it is no longer immersed but has a cusp where its velocity vanishes. This phenomenon was observed by Hill \cite{hill} in his ``moon of maximal lunarity``. Here ``lunarity'' refers to the period of the orbit, which in case of the moon is about the length of a month. Hill used the period as the family parameter which is changing with changing energy.
He observed that for some length of the month which is more than half a year, the velocity of the direct periodic orbit can vanish at some moment. Hill thought that the family of the direct periodic orbit might stop at this moment and wrote
\begin{quote}
Whether this class of satellites is properly to be prolonged beyond this moon, can only be decided by further employment of mechanical quadratures. But it is at least certain that the orbits, if they do exist, do not intersect the line of quadratures, and that the moons describing them would make oscillations to and for, never departing as much as 90 degrees from the point of conjunction or opposition.
\end{quote}
This reasoning was criticized by Adams and Poincar\'e \cite[Chapter 2.3.9]{barrow-green}. In particular, Poincar\'e explained \cite[Chapter III]{poincare} that after increasing the period further the cusps should give rise to little loops. Numerically, the phenomenon predicted by Poincar\'e was verified by Matukuma \cite{matukuma} and H\'enon \cite{henon0}. 

The second disaster which can occur along a family are collisions. However, two-body collisions can be regularized and the family can be continued beyond the collision. In view of these disasters, we see that neither
the rotation number nor the winding number around the primary are preserved in a family of periodic orbits. Nevertheless, in this paper we construct invariants for families of periodic orbits which are not affected by these disasters. 
%The construction of these invariants is based on Arnold's $J^+$-invariant. 

\begin{figure}[!htb]
\centering
\def\svgwidth{1.0\textwidth}%
\begingroup\endlinechar=-1
\resizebox{1.0\textwidth}{!}{%
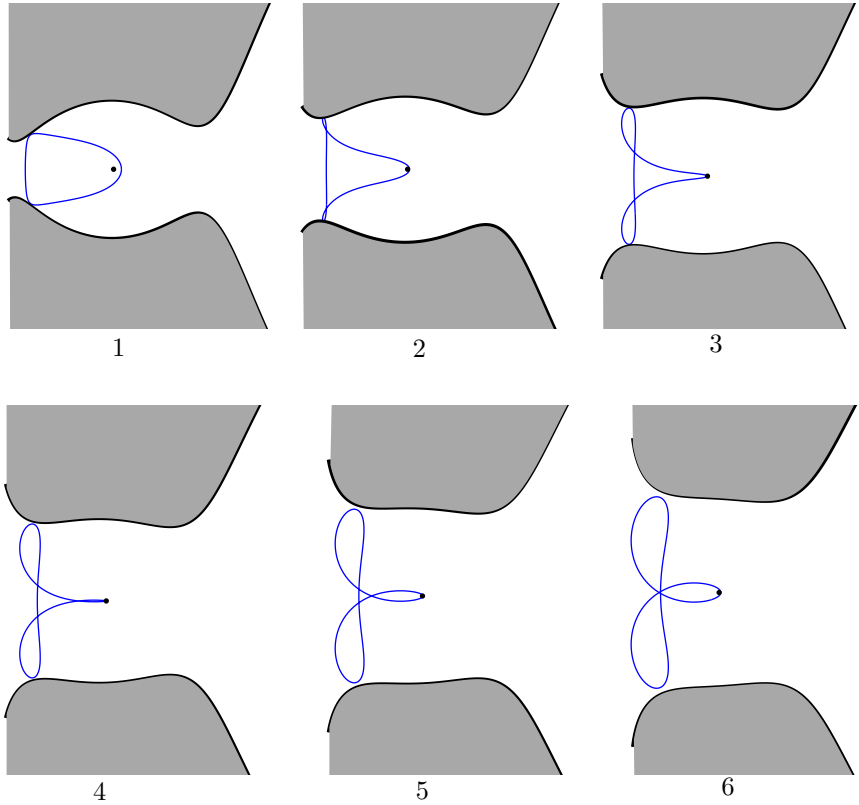%
}\endgroup
\caption{Some members of a $1$-parameter family of periodic orbits in the restricted $3$-body problem ($\mu=0.99$) with Jacobi energy $c$ ranging from 1.575 to 1.55.
When projected to the $q$-plane, the resulting family of immersions goes through several disasters. 
Specifically, between 1 and 2 event $(I_\infty)$ occurs: an exterior loop forms.
Between 3 and 4 event $(I_0)$ occurs: an interior loop forms.
Finally, between 5 and 6 the vertical arc goes through a triple point. This is event $(III)$. These are three of the four possible disasters that can happen in a Stark-Zeeman homotopy, see Definition~\ref{def:Stark-Zeeman-homotopy}.
}
\label{fig:homotopy}
\end{figure}

To conclude this section, we include a family of periodic orbits in the restricted three-body problem, found by numerical integration and a simple shooting argument, that illustrate three of the disasters, see Figure~\ref{fig:homotopy}.

%%%%%%%%%%%%%%%%%%%%%%%%%%%%%%%%%%%%%%%%%%%%%%%%%%%%%%%%%%%%
\section{Planar Stark-Zeeman systems}\label{stark-zeeman}
%%%%%%%%%%%%%%%%%%%%%%%%%%%%%%%%%%%%%%%%%%%%%%%%%%%%%%%%%%%%

%%%
\subsection{Stark-Zeeman systems}
%%%
A {\em Stark-Zeeman system} describes the dynamics of an
electron attracted by a proton subject to an exterior
electric and magnetic field. We suppose that the proton is 
at the origin. Therefore, in dimension $n$ the Coulomb potential of the proton is 
$$V_0 \colon \mathbb{R}^n \setminus \{0\} \to \mathbb{R}, \quad q \mapsto -\frac{1}{|q|}.$$
The exterior electric field is described by a smooth potential
$V_1 \colon U_0 \to \mathbb{R}$, 
where $U_0 \subset \mathbb{R}^n$ is an open subset containing the origin. We abbreviate
$U:=U_0 \setminus \{0\}$. 
The full potential is then the sum
$$V=V_0+V_1 \colon U \to \mathbb{R}.$$
The magnetic field is modeled by a closed $2$-form 
$$
   \sigma_B = \frac12\sum_{i,j=1}^n B_{ij}(q)dq_i\wedge dq_j\in\Omega^2(U_0),
   \qquad B_{ij}=-B_{ji}.
$$ 
It induces the twisted symplectic form
$$\omega_B=\sum_{i=1}^ndp_i \wedge dq_i + \pi^*\sigma_B$$
on the  cotangent bundle $T^*U$, where $\pi \colon T^*U \to U$ is the
footpoint projection. We consider the Hamiltonian
$$H_V \colon T^*U\to\mathbb{R}, \quad (q,p) \mapsto \frac{1}{2}|p|^2+V(q),$$
where the norm is taken with respect to the standard metric on
$\mathbb{R}^n$. 
The Hamiltonian vector field $X^B_V$ is implicitly defined by the condition
$$-dH_V=\omega_B(X^B_V, \cdot).$$
Explicitly, it is given by the formula
$$
   X^B_V = \sum_ip_i\frac{\p}{\p q_i} +
   \sum_i\Bigl(\sum_jB_{ij}(q)p_j-\frac{\p V(q)}{\p q_i}\Bigr)\frac{\p}{\p p_i}.
$$
In particular, in terms of the matrix $B=(B_{ij})$ Hamilton's equation for $X^B_V$ reads
$$
\left\{\begin{array}{l}
\dot q = p\\
\dot p = B(q)p - \nabla V(q),
\end{array}\right.$$
or interpreted as a second order ODE,
\begin{equation}\label{ode}
   \ddot q = B(q)\dot q - \nabla V(q).
\end{equation}
If the summands on the  right hand side are interpreted as forces,
then the potential leads to a force which is only position dependent
but independent of the velocity, while the magnetic field gives 
rise to a force depending linearly on the velocity. 
For $c \in \mathbb{R}$ we abbreviate the energy hypersurface of energy $c$ by
$$\Sigma_c:=H_V^{-1}(c).$$
Then the {\em Hill's region} is defined as
$$\mathfrak{K}_c:=\pi(\Sigma_c)=\{q \in U\mid V(q) \leq c\}.$$
In view of the singularity of $V$ at the origin, the Hill's region has
a unique connected component 
$$\mathfrak{K}_c^b \subset \mathfrak{K}_c$$
which contains $0$ in its closure. We denote by 
$$c_1=c_1(B,V) \in \mathbb{R}\cup\{\infty\}$$
the supremum over all real numbers $c_1$ such for all $c<c_1$ the
following holds: $c$ is a regular value of $V$ (and therefore as well
of $H_V$), and the component $\mathfrak{K}_c^b$ is diffeomorphic to a
closed disk from which an interior point is removed (note that this
always holds for $c$ sufficiently small).  
For example, if $V=V_0$ is the Coulomb potential, then $c_1=0$ is
referred to as the ``escape threshold'' in the physics literature, see for instance \cite{friedrich-wintgen}.
If $V_1$ is defined on the whole space $\mathbb{R}^n$ and goes to
$-\infty$ as $|q|$ goes to infinity, then $c_1$ is the first critical
value of $V$ (respectively $H_V$). 

In this paper we will restrict to the case that $\sigma_B$ is exact
(this is automatic in the planar case $n=2$), so $\sigma_B=dA$ for a
$1$-form 
$$
   A = \sum_jA_j(q)dq_j \in \Omega^1(U_0),\qquad B_{ij}=\frac{\p
     A_j}{\p q_i}-\frac{\p A_i}{\p q_j}. 
$$ 
Then we can untwist the symplectic form by the map
$$\Phi_A \colon T^* U \to T^* U, \quad (q,p) \mapsto
(q,p-A(q)).$$
Indeed, $\Phi_A$ pulls back the twisted symplectic form to the standard symplectic form  
$$\Phi_A^* \omega_{B}=\omega_0=\sum_idp_i\wedge dq_i$$
and the Hamiltonian $H_V$ to 
$$
   H^A_V \colon T^*U\to\mathbb{R}, \quad (q,p) \mapsto \frac{1}{2}|p-A(q)|^2+V(q),   
$$
so Hamilton's equations of $H_V$ with respect to the twisted
symplectic form $\omega_B$ transform into Hamilton's equations of the
twisted Hamiltonian $H_V^A$ with respect to the untwisted form $\omega_0$. Explicitly, 
the Hamiltonian vector field $X_V^A$ of the Hamiltonian 
with respect to the untwisted symplectic form $\omega_0$ equals
$$
   X^A_V = \sum_i\Bigl(p_i-A_i(q)\Bigr)\frac{\p}{\p q_i} +
   \sum_i\Bigl(\sum_j\frac{\p A_j(q)}{\p q_i}(p_j-A_j(q))-\frac{\p V(q)}{\p q_i}\Bigr)\frac{\p}{\p p_i},
$$
hence Hamilton's equation for $X^A_V$ reads
$$
\left\{\begin{array}{l}
\dot q_i = p_i-A_i(q)\\
\dot p_i = \sum_j\frac{\p A_j(q)}{\p q_i}(p_j-A_j(q)) - \frac{\p
  V(q)}{\p q_i},
\end{array}\right.$$
which gives rise to the same second order ODE~\eqref{ode}.

{\bf The planar case. }
A {\em planar Stark-Zeeman system} is a Stark-Zeeman system with $n=2$
degrees of freedom. In this case the magnetic field can be written as
$$
   \sigma_B = B(q)dq_2\wedge dq_1
$$
for a smooth function $B:U_0\to\R$ 
(so our convention is such that
positive $B$ corresponds to a magnetic field pointing in the negative
$q_3$-direction). The corresponding $2\times 2$-matrix is given by $B(q)J$, where 
$$J=\left(\begin{array}{cc}
0 & -1\\
1 & 0
\end{array}\right)$$
is the matrix of complex multiplication by $i$ if $\mathbb{R}^2$ is
identified with $\mathbb{C}$. Thus Newton's equation~\eqref{ode}
reads 
\begin{equation}\label{ode2}
   \ddot q = BJ \dot q -\nabla V(q),
\end{equation}
or written out in components,
$$
\left\{\begin{array}{l}
   \ddot q_1 = -B \dot q_2 - \tfrac{\partial V}{\partial q_1}\\
   \ddot q_2 = B \dot q_1 -\frac{\partial V}{\partial q_2}.
\end{array}\right.
$$

%%%
\subsection{Examples of planar Stark-Zeeman systems}\label{examples}
%%%

Planar Stark-Zeeman systems include many systems of physical relevance
in classical mechanics as well as in quantum mechanics. We consider here
only the classical systems. In all the examples the magnetic field $B$
will be constant. The integrability of some of these systems is discussed in Section~\ref{sec:sep} below.

{\bf Kepler problem. }
The case $B=V_1=0$ corresponds to the Kepler problem of two bodies
moving under their mutual gravitation, or equivalently, a (classical)
proton and electron moving under their mutual electric attraction. 

{\bf Stark problem. }
The case $B=0$, $V_1(q)=-Eq_1$ corresponds to the $2$-dimensional Stark problem of an
electron moving in the field of a proton and a constant electric field
of strength $E$ in $q_1$-direction. This system is integrable.

{\bf Zeeman problem. }
%\marginpar{Does this make physical sense?}
The case $V_1=0$ and constant $B\neq 0$ corresponds to the $2$-dimensional
Zeeman problem of an electron in the plane moving in the field of a
proton and a constant magnetic field perpendicular to the plane. 

{\bf Diamagnetic Kepler problem. }
Consider the $3$-dimensional Zeeman problem of an of an electron in
$\R^3$ with coordinates ${\bf r}=(x,y,z)$ moving in the field of a
proton at the origin and a constant magnetic field $B\p_z$ in
$z$-direction. This field has the vector potential $A =
\frac{B}{2}(x\p_y-y\p_x)$, so the twisted Hamiltonian is 
\begin{align*}
   H_V^A &= \frac12|p-A|^2 - \frac{1}{|{\bf r}|} \cr
   &= \frac{|p|^2}{2} - \frac{B}{2}(xp_y-yp_x) +
   \frac{B^2}{8}(x^2+y^2) - \frac{1}{|{\bf r}|} \cr
   &= \frac{|p|^2}{2} - \frac{B}{2}(xp_y-yp_x) +
   \frac{B^2}{8}(x^2+y^2) - \frac{1}{|{\bf r}|} \cr
   &= \frac12\Bigl(p_\rho^2+p_z^2+\frac{p_\theta^2}{\rho^2} -
   Bp_\theta\Bigr) + \frac{B^2}{8}\rho^2 - \frac{1}{\sqrt{\rho^2+z^2}},
\end{align*}
where $(\rho,\theta,z)$ are cylindrical coordinates and
$p_\theta=xp_y-yp_x$ is the $z$-component of angular momentum. 
Note that $p_\theta$ is a conserved quantity. Restricting to the case
$p_\theta=0$ and setting $q_1=\rho$ and $q_2=z$ we obtain a planar
Stark-Zeeman system with vanishing magnetic field and
$V_1=\frac{B^2}{8}q_1^2$, i.e.~with twisted Hamiltonian 
$$
   H_V^A = \frac12|p|^2-\frac{1}{|q|}+\frac{B^2}{8}q_1^2. 
$$
While the preceding systems were variations of the classical hydrogen
atom, the following systems will describe variations of the $2$-body
problem in celestial mechanics. 

{\bf Rotating Kepler problem. }
Consider a coordinate system $(q_1,q_2)$ that is
rotating with angular velocity $\om$ in counterclockwise direction with
respect to the fixed coordinate system $(q_1',q_2')$. Thus in complex
notation $q=q_1+iq_2$ and $q'=q_1'+iq_2'$ and their time derivatives
are related by  
\begin{gather*}
   q=e^{i\om t}q',\qquad \dot q=e^{i\om t}(\dot q'+i\om q'), \cr
   \ddot q = e^{i\om t}(\ddot q'+2i\om\dot q'-\om^2q') 
   = e^{i\om t}F_t(e^{-i\om t}q)  + 2i\om\dot q + \om^2q,
\end{gather*}
where $F_t$ is the (possibly time-dependent) force in Newton's law
$\ddot q'=F_t(q')$. In the Kepler problem, $F(q)=-q/|q|^3=e^{i\om
  t}F_t(e^{-i\om t}q)$ and Newton's law in rotating coordinates reads
$$
   \ddot q = -\frac{q}{|q|^3} + 2i\om\dot q + \om^2q,
$$
which is a planar Stark-Zeeman system with magnetic field $2\om$
(generating the Coriolis force) and potential
$V_1=-\frac12\om|q|^2$ (generating the centrifugal force). 

{\bf Restricted $3$-body problem. }
The planar circular restricted $3$-body problem describes the motion
$q(t)$ of a massless body in the plane in the (time-dependent)
gravitational field generated by two massive bodies (say, the sun and
the earth of masses $M_S$ and $M_E$) moving on circular orbits around
their center of mass. We put the center of mass at the origin and
write the motion of the sun and earth using complex notation as 
$$
   q_S(t)=e^{-i\om t}q_S,\qquad q_E(t)=e^{-i\om t}q_E
$$
with fixed $q_S,q_E\in\C$. Newton's law for their Keplerian motion yields
the relations
$$
   \om^2 = \frac{M_S}{|q_E|\,|q_S-q_E|^2} = \frac{M_E}{|q_S|\,|q_S-q_E|^2}.   
$$
In particular, $M_S/M_E = |q_E|/|q_S|$. So we can uniquely write
$$
   M_E=\mu M,\quad M_S=(1-\mu)M,\quad q_E=-R(1-\mu),\quad q_S=R\mu
$$
for the total mass $M=M_S+M_E$ and parameters $\mu,R>0$, were we have
chosen $q_E$ on the negative real axis and $q_S$ on the positive real axis.
The relation for $\om$ then becomes
\begin{equation}\label{eq:Kepler3}
   \om^2 = \frac{M}{R^3},
\end{equation}
which just expresses Kepler's 3rd law. The motions of the massive
bodies read 
$$
   q_E(t)=-e^{-i\om t}R(1-\mu),\qquad q_S(t)=e^{-i\om t}R\mu
$$
and generates the gravitational force 
$$
  F_t(q) = -\frac{\mu M(q-q_E(t))}{|q-q_E(t)|^3} -\frac{(1-\mu)M(q-q_S(t))}{|q-q_S(t)|^3}
$$
on the massless body. We switch to rotating coordinates as in the
rotating Kepler problem. The gravitational force then becomes time-independent,
$$
   e^{i\om t}F_t(e^{-i\om t}q) = -\frac{\mu M(q+R(1-\mu))}{|q+R(1-\mu)|^3} -
   \frac{(1-\mu)M(q-R\mu)}{|q-R\mu|^3},
$$
and adding in the Coriolis and centrifugal forces we obtain Newton's
equation in rotating coordinates,
$$
   \ddot q = 2\om i\dot q + \om^2q -\frac{\mu M(q+R(1-\mu))}{|q+R(1-\mu)|^3} -
   \frac{(1-\mu)M(q-R\mu)}{|q-R\mu|^3}.
$$
This equation is generated by the magnetic field $B=2\om$ and the
potential
$$
   V(q) = -\frac12\om^2|q|^2 -\frac{\mu M}{|q+R(1-\mu)|} -
   \frac{(1-\mu)M}{|q-R\mu|},
$$
where $\om$ is given by~\eqref{eq:Kepler3}. Normalizing $M=R=\om=1$
this yields the familiar description of the planar circular restricted
$3$-body problem. This becomes a Stark-Zeeman system if we shift the
origin to one of the massive bodies, say the sun. Replacing $q$ by
$q-R(1-\mu)$ then yields the equivalent potential
$$
   \wt V(q)=V(q-R(1-\mu)) = -\frac12\om^2|q-R(1-\mu)|^2 -\frac{\mu M}{|q|} -
   \frac{(1-\mu)M}{|q-R|}. 
$$
{\bf Hill's lunar problem. }
Hill's lunar problem is the limiting case of the planar circular
restricted $3$-body problem as $M,R\to\infty$ and $\mu\to 0$ subject
to the restrictions
\begin{equation}\label{eq:Hill-limit}
   \frac{M}{R^3}=\om^2\to {\rm constant}\neq 0,\qquad \mu M\to 1.  
\end{equation}
Using the Taylor expansion
$$
   \frac{1}{\sqrt{1+x}} = 1-\frac{1}{2}x + \frac{3}{8}x^2 + O(x^3)
$$ 
we expand the gravitational potential of the sun in $\wt V$ to second order in $1/R$:
\begin{align*}
   \frac{(1-\mu)M}{|q-R|} 
   &= \frac{(1-\mu)M}{R}\frac{1}{\sqrt{1-\frac{2q_1}{R}+\frac{q^2}{R^2}}} \cr
   &= \frac{(1-\mu)M}{R}\Bigl(1-\frac{|q|^2}{2R^2}+\frac{q_1}{R}+\frac{3}{2}\frac{q_1^2}{R^2}
   +O(\frac{1}{R^3})\Bigr).
\end{align*}
The conditions~\eqref{eq:Hill-limit} are chosen such that the terms
with $q_i$ and $|q|^2$ cancel with the corresponding terms in the
centrifugal potential
$$
   \frac12\om^2|q-R(1-\mu)|^2 = \frac{M}{2R^3}\Bigl(|q|^2-2R(1-\mu)q_1+R^2(1-\mu)^2\Bigr),
$$
so dropping constant terms we obtain in the limit the potential
$$
   V_{\rm Hill}(q) = -\frac{1}{|q|}-\frac{3}{2}\om^2q_1^2. 
$$
Together with the magnetic field of strength $B=2\om$ this describes
Hill's lunar problem. By constant rescaling in time we can normalize
$\om$ to $1$. 

{\bf Frozen Hill's problem. }
Dropping the magnetic field in Hill's problem gives rise to the
``frozen Hill's problem''. The frozen Hill's problem can be interpreted
as a diamagnetic Kepler problem for an imaginary magnetic field $B$.

{\bf Frozen Hill's problem with centrifugal force. }
Surprisingly, adding a suitable centrifugal term to the frozen Hill's
problem leads to an integrable system described by the Hamiltonian
$$H^A_V=\tfrac{1}{2}|p|^2-\tfrac{1}{|q|}-\tfrac{3}{2}q_1^2-\tfrac{1}{2}|q|^2.$$

{\bf Euler problem. }
Dropping the magnetic field and the centrifugal potential from the
restricted $3$-body problem leads to the Euler problem with Hamiltonian
$$
   H^A_V = \frac12|p|^2 -\frac{\mu M}{|q+R(1-\mu)|} -
   \frac{(1-\mu)M}{|q-R\mu|}.
$$
It describes the motion of an electron in the electric field of two
protons held at fixed position. This system is integrable. 

{\bf Lagrange problem. }
Adding an additional elastic force with center between the two masses to the potential of the Euler problem leads to
another integrable system first considered by Lagrange. Namely, for any $a \in \mathbb{R}$ the 
Hamiltonian 
$$
   H^A_V = \frac12|p|^2 -\frac{\mu M}{|q+\frac{1}{2}R|} -
   \frac{(1-\mu)M}{|q-\frac{1}{2}R|}+a|q|^2.
$$
is integrable. 
If we think of this force as a centrifugal force, then in the case where the two masses are equal this corresponds to the restricted $3$-body problem with Coriolis force removed. Note
that in this case the Lagrange problem and the restricted $3$-body problem have the same potential and therefore the same Hill's regions. 
%To do: 
%Discuss integrability / chaos of the systems (Waldvogel, Diss
%Winterberger for Hill, Xia, Poincare for restricted 3BP\\
%Explain the limits Euler $\to$ Stark and frozen Hill with centrifugal
%$\to$ Lagrange, and their integrability (Hiertarinta).\\
%What is consequence of integrability to second order in diamagnetic
%Kepler for KAM theory?

\subsection{Integrability}\label{sec:sep}

In this section we explain how the integrability of various systems in the preceding section can be derived in a unified manner from results in \cite{landau-lifschitz,hietarinta}.
Consider a mechanical Hamiltonian, i.e., a Hamiltonian just consisting of kinetic and potential energy,
$$H(q,p)=\frac{1}{2}p^2+V(q)$$
where $V \colon U \to \mathbb{R}$ is the potential for $U \subset \mathbb{R}^2$ an open subset. 
Abbreviate $r=|q|$. The following two well-known theorems gives some conditions on the potential $V$ which ensure that the Hamiltonian $H$ admits an integral which is quadratic in the momenta. 
The following two theorems can be found for example in \cite[Section 48]{landau-lifschitz} or \cite{hietarinta}.
\begin{thm}\label{sep1}
Assume that the potential $V$ can be written in the form
$$V=\frac{f(r+q_1)+g(r-q_1)}{r}$$
for smooth functions $f, g \colon \mathbb{R} \to \mathbb{R}$. Then the function
$$I=(q_1p_2-q_2p_1)p_2+\frac{(r+q_1)g(r-q_1)-(r-q_1)f(r+q_1)}{r}$$
is an integral for the Hamiltonian $H$, i.e., $H$ and $I$ Poisson commute. 
\end{thm}
There are other conditions on the potential $V$ which guarantee an integral quadratic in the momenta.
For $c>0$ we abbreviate
$$r_1=\sqrt{(q_1+c)^2+q_2^2}, \quad r_2=\sqrt{(q_1-c)^2+q_2^2}$$
and further set
$$u=\tfrac{1}{2}(r_1+r_2), \quad v=\tfrac{1}{2}(r_1-r_2).$$
%We are now in position to state
\begin{thm}\label{sep2}
Assume that the potential $V$ can be written in the form
$$V=\frac{f(u)-g(v)}{u^2-v^2}$$
for smooth functions $f,g \colon \mathbb{R} \to \mathbb{R}$. Then the function
$$I=(q_1p_2-q_2p_1)^2+c^2p_1^2+2\frac{v^2f(u)-u^2g(v)}{u^2-v^2}$$
is an integral for the Hamiltonian $H$.
\end{thm}
By an elementary but tedious computation one can check directly that under the assumptions of
Theorem~\ref{sep1} and Theorem~\ref{sep2} one has $\{H,I\}=0$. The reader might wonder
how people could have come up with these conditions and found the integral $I$. We refer to \cite{hietarinta}
for an enlightening discussion concerning these questions.  Alternatively Hamilton-Jacobi theory leads
to a proof of the two theorems above as explained in \cite{landau-lifschitz}.

Using the above two theorems we can now readily check the integrability of several of the examples of
Stark-Zeeman systems discussed in Section~\ref{examples}.

{\bf Stark problem. }
Plug $f(s)=-\tfrac{1}{2}-\tfrac{E}{4}s^2$ and $g(s)=-\tfrac{1}{2}+\tfrac{E}{4}s^2$ into 
Theorem~\ref{sep1} to obtain $V=-\tfrac{1}{r}-Eq_1$.

{\bf Frozen Hill's problem with centrifugal force. }
Plug $f(s)=-\tfrac{1}{2}-\tfrac{1}{4}s^3$ and $g(s)=-\tfrac{1}{2}+s^3$ into Theorem~\ref{sep1} to obtain
$V=-\tfrac{1}{r}-\tfrac{3}{2}q_1^2-\tfrac{1}{2}r^2$.

{\bf Euler problem. }
For $c=\tfrac{R}{2}$ plug $f(s)=-Ms$ and $g(s)=(1-2\mu)Ms$ into Theorem~\ref{sep2}  to obtain
$V=-\tfrac{\mu M}{r_1}-\tfrac{(1-\mu)M}{r_2}$, which is the Euler problem up to a translation of the position coordinates.

{\bf Lagrange problem. }
Concerning the Lagrange problem we compute
$$\frac{u^4-v^4}{u^2-v^2}=u^2+v^2=\frac{1}{4}(r_1^2+r_2^2)=\frac{1}{2}(r^2+c^2).$$
In view of Theorem~\ref{sep2}, this means that we can add to the Euler problem an elastic force with center in the middle of the two masses and still keep the system integrable.

\subsection{Levi-Civita regularization of a planar Stark-Zeeman system}\label{levicivita}

We consider a planar Stark-Zeeman system 
$$H=H_V^A \colon T^* U \to \mathbb{R}, \quad (q,p) \mapsto \frac{1}{2}|p-A(q)|^2
+V(q)$$
where $U=U_0 \setminus \{0\}$ for $U_0 \subset \mathbb{C}$ an open subset containing the origin,
$$V=V_0+V_1 \colon U \to \mathbb{R}$$
where $V_0$ is the Coulomb potential and $V_1$ is a smooth function on $U$ extending smoothly to $U_0$, and
$$A \in \Omega^1(U_0)$$
such that $dA=\sigma_B \in \Omega^2(U)$ is the magnetic field. 

Abbreviate $\mathbb{C}^*=\mathbb{C} \setminus \{0\}$. The \emph{Levi-Civita map} is the
transformation
$$\mathcal{L} \colon T^* \mathbb{C}^* \to T^* \mathbb{C}^*$$
given by
$$\mathcal{L}(v,u)=\big(v^2, u/2\bar{v}\big), \quad
(v,u) \in \mathbb{C}^* \times \mathbb{C}=T^* \mathbb{C}^*.$$
It is the cotangent lift of the complex squaring map
$$L \colon \mathbb{C} \to \mathbb{C}, \quad v \mapsto v^2,$$
In particular, the Levi-Civita map pulls back the canonical symplectic structure of $T^* \mathbb{C}$ to itself, i.e., with $(q,p)=\big(v^2, u/2\bar{v}\big)$ we have
$$\mathcal{L}^*(dq_1 \wedge dp_1+dq_2 \wedge dp_2)=dv_1 \wedge du_1+dv_2 \wedge du_2.$$
Set
$$\mathcal{U}:=L^{-1}(U) \subset \mathbb{C}^*.$$
For $c \in \mathbb{R}$ we introduce the Hamiltonian
$$\mathcal{H}_c \colon T^* \mathcal{U} \to \mathbb{R}, \quad
(v,u) \mapsto 4|v|^2\big(H \circ \mathcal{L}(v,u)-c\big).$$
Since the Levi-Civita map is a symplectomorphism, 
it maps up to reparametrization flow lines of the Hamiltonian flow of $\mathcal{H}_c$ of energy
zero to flow lines of the Hamiltonian flow of $H$ of energy $c$. 

Our first goal is to understand the structure of the Hamiltonian $\mathcal{H}_c$. For that purpose we introduce 
some notation. 
First observe that $0$ lies in the closure of $\mathcal{U}$ and we abbreviate
$$\mathcal{U}_0=\mathcal{U} \cup \{0\}.$$
Note that $\mathcal{U}_0$ is an open subset of $\mathbb{C}$.

Using the canonical identification of the tangent and cotangent bundle of $\mathbb{C}$ we can think
of the one-form $A$ as well as a vector field $A \colon U_0 \to \mathbb{C}$ and we introduce 
the one-form
$$\mathcal{A} \in \Omega^1(\mathcal{U}_0)$$
as the vector field
$$\mathcal{A} \colon \mathcal{U}_0 \to \mathbb{C}, \quad
v \mapsto 4\bar{v} A(v^2).$$
We further introduce the potential 
$$\mathcal{V}_c \colon \mathcal{U}_0 \to \mathbb{R}, \quad
v \mapsto -4c|v|^2+4|v|^2 V_1(v^2)-4$$
where we recall that $V_1=V-V_0$ where $V_0$ is the Coulomb potential. 

With these preparations we compute for $(v,u) \in T^* \mathcal{U}$
\begin{eqnarray}\label{regham}
\mathcal{H}_c(v,u)&=&\frac{4|v|^2}{2}\bigg|\frac{u}{2\bar{v}}-A(v^2)\bigg|^2+4|v|^2V(v^2)-4c|v|^2\\ \nonumber
&=&\frac{1}{2}\big|u-4\bar{v}A(v^2)\big|^2-4c|v|^2+4|v|^2 V_1(v^2)-4\\ \nonumber
&=&\frac{1}{2}\big|u-\mathcal{A}(v)\big|^2+\mathcal{V}_c(v).
\end{eqnarray}
We first observe that, although we just defined $\mathcal{H}_c$ on $T^* \mathcal{U}$, both
maps $\mathcal{A}$ and $\mathcal{V}_c$ are actually defined smoothly on $\mathcal{U}_0$, so that 
we can smoothly extend $\mathcal{H}_c$ to $T^* \mathcal{U}_0$. In the following we think by abuse of notation of
$$\mathcal{H}_c \colon T^* \mathcal{U}_0 \to \mathbb{R}$$
as the smooth map given by formula (\ref{regham}) and refer to it as the \emph{Levi-Civita regularized Hamiltonian}. The fibre over the origin has the following physical significance: it corresponds to collisions
of the particle with the charge at the origin. The Hamiltonian $\mathcal{H}_c$ regularizes these collisions after
time reparametrization, so that the regularized orbits just pass through the origin.

Another interesting observation is that the Levi-Civita regularized Hamiltonian is still a magnetic Hamiltonian, consisting of a twisted kinetic energy term and a potential term. In contrast to $H$, the potential in $\mathcal{H}_c$ has no singularity at the origin. 

Abbreviate the energy hypersurface of the Levi-Civita regularized Hamiltonian by
$$\mathcal{S}_c=\mathcal{H}_c^{-1}(0).$$
As for the nonregularized Hamiltonian we define the 
\emph{Hill's region} for the regularized Hamiltonian as
$$\mathfrak{K}_c^{\mathrm{reg}}=\pi(\mathcal{S}_c)=
\{v \in \mathbb{C}: \mathcal{V}_c(v) \leq 0\} \subset \mathcal{U}_0$$
where $\pi \colon T^* \mathcal{U}_0 \to \mathcal{U}_0$ is the footpoint projection. Note that the regularized Hill's region is related to the nonregularized Hill's region by
$$\mathfrak{K}_c^{\mathrm{reg}}=L^{-1}(\mathfrak{K}_c) \cup \{0\}.$$
Recall that $c_1 \in \mathbb{R} \cup \{\infty\}$ denotes the first critical value respectively the escape threshold. If $c<c_1$ the Hill's region contains a unique bounded component 
$\mathfrak{K}_c^b \subset \mathfrak{K}_c$ which is characterized by the property that the origin lies in its closure. Moreover, $\mathfrak{K}_c^b$ is diffeomorphic to a pointed closed disk. We define
$$\mathfrak{K}_c^{b, \mathrm{reg}} \subset \mathfrak{K}_c^{\mathrm{reg}}$$
as the connected component of $\mathfrak{K}_c^{\mathrm{reg}}$ which contains the origin. If $c<c_1$
the two components are related to each other by
$$\mathfrak{K}_c^{b,\mathrm{reg}}=L^{-1}(\mathfrak{K}_c^b) \cup \{0\}.$$
In particular, $\mathfrak{K}_c^{b,\mathrm{reg}}$ is diffeomorphic to a closed disk. 

We denote by
$$\mathcal{S}_c^b=\{(v,u) \in \mathcal{S}_c: v \in \mathfrak{K}_c^{b,\mathrm{reg}}\}$$
the connected component of the energy hypersurface $\mathcal{S}_c^b$ which projects to 
$\mathfrak{K}_c^{b,\mathrm{reg}}$. Note that for each point $v$ in the interior 
of the Hill's region $\mathfrak{K}_c^{b,\mathrm{reg}}$ the fiber in $\mathcal{S}_c^b$ lying over $v$ 
is a circle, where at the zero velocity curves (i.e., the boundary of 
$\mathfrak{K}_c^{b,\mathrm{reg}}$) the fiber collapses to a point. We deduce that topologically
$\mathcal{S}_c^b$ is a $3$-dimensional sphere. 

We summarize the preceding discussion as follows.

\begin{prop}
The Levi-Civita regularized Hamiltonian at energy $c<c_1$ is a nonsingular magnetic Hamiltonian $\mathcal{H}_c:T^* \mathcal{U}_0 \to \mathbb{R}$. The bounded component $\mathfrak{K}_c^{b,\mathrm{reg}}$ of its Hill's region is a closed disk whose preimage $\mathcal{S}_c^b\subset \mathcal{H}_c^{-1}(0)$ is a $3$-sphere. \hfill$\square$
\end{prop}

%%%
\subsection{Moser regularization of a Stark-Zeeman system}
%%%
Let us consider a Stark-Zeeman system in any dimension $n$. 
Even if we restrict our attention for an energy value $c<c_1$ to the
connected component of the energy hypersurface above the bounded
component $\mathfrak{K}_c^b$ of the Hill's region, 
$$\Sigma_c^b:=\{(q,p) \in \Sigma_c\mid q \in \mathfrak{K}_c^b\},$$
we are facing noncompactness issues because the electron might
collide with the proton. However, two-body collisions can always be
regularized as we now recall.  
For any regular energy value $c$ of $H_V$ we introduce the Hamiltonian
$K_c \colon T^* U \to \mathbb{R}$ given by
\begin{align*}
   K_c(q,p) &:= |q|\big(H_V^A(q,p)-c\big) \cr
   &= \frac{|q|}{2}|p-A(q)|^2 - 1 + |q|V_1(q) - c|q| \cr 
   &= \frac{|q|}{2}|p|^2 + |q|\Bigl(\la p,A(q)\ra+\frac12|A(q)|^2 + V_1(q) -
   c\Bigr) - 1 
\end{align*}
Since $q\neq 0$ on $U$, the level sets satisfy 
\begin{equation}\label{regu}
   K_c^{-1}(0) = (H_V^A)^{-1}(c) = \Phi_A^{-1}(\Sigma_c).
\end{equation}
Note that $|p|\to\infty$ as $q\to 0$ in this level set. To study the
behaviour as $p$ approaches infinity, we pull back $K_c$ under the simple symplectic map
$$\sigma \colon T^* \mathbb{R}^n \to T^* \mathbb{R}^n, \quad (q,p) \mapsto (-p,q).$$
Note that this map is from the physical point of view rather
revolutionary because it interchanges the roles of base and fiber, 
so positions become momenta and vice versa. We will keep the original
notation $(q,p)$ and just remember that the $p_i$ are now interpreted
as positions and the $q_i$ as negative momenta. 

To study the behaviour as $|p|\to\infty$, we view $\mathbb{R}^n$ as a
chart of the sphere $S^n_r\subset\R^{n+1}$ of radius $r>0$ via stereographic
projection 
$$
   \psi_r^N:\R^n\stackrel{\cong}\longrightarrow S^n_r\setminus\{N\},
   \quad p\mapsto(x,x_{n+1}), 
$$
where $N=(0,\dots,0,r)$ is the north pole and $x=(x_1,\dots,x_n)\in\R^n$. 
By elementary geometry, the map $\psi_r^N$ is given by
\begin{equation}\label{eq:stereographic}
   x = \frac{2r^2p}{|p|^2+r^2},\qquad x_{n+1} =
   r\frac{|p|^2-r^2}{|p|^2+r^2}. 
\end{equation}
To see whether the Hamiltonian $K_c$ smoothly extends to the fiber at infinity of
$T^*S^n$ we consider the transition map 
$$
   \psi_r:=(\psi_r^S)^{-1}\circ\psi_r^N \colon \mathbb{R}^n \setminus \{0\} \to \mathbb{R}^n \setminus \{0\}, \quad
   p \mapsto \frac{r^2p}{|p|^2}
$$
between the stereographic projections from the north and south poles. 
%Let us for now restrict to the case $r=1$ and abbreviate $\psi=\psi_1$. 
The Jacobi matrix $D\psi_r(p)$ has entries
$$
   D\psi_r(p)_{ij} = \frac{r^2}{|p|^4}\begin{cases}
   -2p_ip_j & i\neq j \\
   |p|^2-2p_i^2 & i= j. 
   \end{cases}
$$
Since this matrix is symmetric and $\psi_r\circ\psi_r=\id$, its
inverse transpose computes to be (using $|\psi_r(p)|=r^2/|p|$)
$$
   \big(D\psi_r(p)^{-1}\big)_{ij}^T = D\psi_r\bigl(\phi_r(p)\bigr)_{ij} 
   = \frac{1}{r^2}\begin{cases}
   -2p_ip_j & i\neq j \\
   |p|^2-2p_i^2 & i= j. 
   \end{cases}
$$
Therefore, the induced symplectomorphism
$$
   \Psi_r \colon T^*(\mathbb{R}^n \setminus \{0\}) \to T^*(\mathbb{R}^n
   \setminus \{0\}),\quad (p,q)\mapsto
   \Bigl(\psi_r(p),\big(D\psi_r(p)^{-1}\big)_{ij}^T\cdot q\Bigr) =:(\tilde p,\tilde q)
$$
is given by
$$
   \tilde p=\frac{r^2p}{|p|^2},\qquad \tilde q=\frac{|p|^2q-2\la p,q\ra p}{r^2}.
$$
In particular $|\tilde q|=|p|^2|q|/r^2$, and using this we obtain for $(p,q) \in \Psi_r^{-1}(T^*U)$
$$
   K_c\circ\Psi_r(p,q) = \frac{r^2}{2}|q| + |q|\la p,A(\tilde q)\ra
   + \frac{|p|^2|q|}{r^2}\Bigl(\frac{|A(|\tilde q)|^2}{2} + V_1(\tilde
   q)-c\Bigr)-1. 
$$
We see that this extends continuously over $p=0$ via
$$
   K_c\circ\Psi_r(0,q) = \frac{r^2}{2}|q|-1.
$$
In particular, on the level set $K_c^{-1}(0)$ we have $|q|\to 2/r^2$ as
$p\to 0$, so the extension of $K_c\circ\Psi_r$ is smooth near this
level set. Hence the hypersurface
$$
   \Phi_r^{-1}\bigl(K_c^{-1}(0)\bigr) =
   \Psi_r^{-1}\circ\Phi_A^{-1}(\Sigma_c) \subset T^*(\R^n\setminus\{0\})
$$
can be extended smoothly over $p=0$ by adding the sphere
$\{p=0,\,|q|=2/r^2\}$. 
Now a straightforward computation using~\eqref{eq:stereographic} shows
that the pullback of the round metric $g_r$ on the sphere
$S^n_r\subset\R^{n+1}$ of radius $r$ under the stereographic
projection $\psi_r^N:\R^n\to S^n_r\setminus\{N\}$ is given at the
point $p\in\R^n$ by 
$$
   \bigl((\psi_r^N)^*g_r\bigr)_p = \frac{4r^4}{(|p|^2+r^2)^2}\la\ ,\ \ra. 
$$
So the sphere $\{p=0,\,|q|=2/r^2\}\subset T_0^*\R^n$ is the preimage
under $\Psi_r^N$ of the sphere in $T_N^*S^n_r$ of radius $1/r^2$ with
respect to the dual of the round metric $g_r$. We summarize this
discussion in

\begin{prop}\label{prop:Moser}
For any regular value $c$ of $H_V$ and any $r>0$ the closure of the hypersurface
$\Psi_r^N\circ\Phi_A^{-1}(\Sigma_c)\subset T^*S^n_r$ is a smooth
hypersurface $S_{c,r}\subset T^*S^n_r$ which intersects the
fibre over the north pole in the sphere of radius $1/r^2$ with
respect to the dual of the round metric $g_r$ on $S^n_r$. \hfill$\square$
\end{prop}

\begin{ex}\label{ex:Kepler}
{\rm An illustrative example is the case where there is neither an exterior
electric nor magnetic field, i.e., just the case of an electron
attracted by a proton which corresponds to the Kepler problem. Then 
$$
   H(q,p)=\frac12|p|^2-\frac{1}{|q|},\qquad K_c = |q|\Bigl(\frac12|p|^2-c\Bigr)-1.
$$
Suppose $c<0$ and set $r:=\sqrt{-2c}$. Then the energy hypersurface
$$
   \Sigma_c=K_c^{-1}(0) = \Bigl\{(p,q)\;\Bigl|\; \frac{|p|^2+r^2}{2r^2}|q|=\frac{1}{r^2}\Bigr\}
$$
corresponds under stereographic projection to the radius $1/r^2$
sphere bundle in $T^*S^n_r$ (where we interpret $p$ as position and
$q$ as momentum in $T^*S^n_r$). So the regularized energy hypersurface
$\mathcal{S}_{c,r}$ is the radius $1/r^2$ sphere bundle in $T^*S^n_r$
and the regularized Hamiltonian
flow just becomes (up to time reparametrization) the geodesic flow for
the round metric on the sphere $S^n_r$ of radius $r=\sqrt{-2c}$. In
particular, for $c=-1/2$ we get the geodesic flow on the unit sphere.}
\end{ex}

The preceding example was the main reason for allowing for an
arbitrary radius $r>0$ in our discussion. From now on we will restrict
to the case $r=1$ and drop $r$ from the notation. 

In general, suppose that $c<c_1$ and consider the connected component
$\Sigma_c^b$ of the energy hypersurface lying above the bounded part
$\mathfrak{K}_c^b$ of the Hill's region which contains the origin in its
closure. This component is compact up to collisions of the electron
with the proton, which after regularization correspond to points in
the cotangent sphere over the north pole in $S^n$. 
Hence the corresponding connected component
$$\mathcal{S}_{c}^b \subset \mathcal{S}_{c}$$
of the regularized hypersurface $\mathcal{S}_{c}\subset T^*S^n$ is
{\em compact}. We can homotop $\mathcal{S}_c^b$ through regular compact
energy hypersurfaces to Example~\ref{ex:Kepler} by first decreasing
$c$ to a very negative value and then switching off the
external electric and magnetic field. In particular, the energy
hypersurfaces in such a homotopy remain diffeomorphic so that we
proved  

\begin{cor}\label{cor:Moser}
If $c<c_1$, then the bounded component $\mathcal{S}_c^b$ of the Moser regularized
energy hypersurface is diffeomorphic to the unit cotangent bundle
$S^*S^n$. \hfill$\square$
\end{cor}

\begin{rem}
The Levi-Civita regularized energy hypersurface is the preimage of the Moser regularized energy hypersurface under a suitable $2$--$1$ cover $S^3\to\R P^3\cong S^*S^2$. 
\end{rem}

%%%%%%%%%%%%%%%%%%%%%%%%%%%%%%%%%%%%%%%%%%%%%%%%%%%%%%%%%%%%
\section{Topology of periodic orbits}\label{topology}
%%%%%%%%%%%%%%%%%%%%%%%%%%%%%%%%%%%%%%%%%%%%%%%%%%%%%%%%%%%%

%Let $U \subset \mathbb{R}^2$ be an open subset. We consider the
%Hamiltonian system with magnetic field given by the smooth function
%$B:U\to\mathbb{R}$ and a smooth potential $V:U\to\mathbb{R}$. In this
%section, we need not assume that the system is of Stark-Zeeman type. 

%%%
\subsection{Periodic orbits in planar Stark-Zeeman systems}
%%%

{\bf Local structure of orbits. }
Consider an orbit $q:\R\to\R^2$ of a planar Stark-Zeeman system of energy $c<c_1$. We allow $q$ to pass through the origin, where it is understood to be the projection of a regularized collision orbit. Let us describe the structure of $q$ near a point $q_0=q(t_0)$. We distinguish $3$ cases.

{\em Case 1: }$q_0\neq 0$ lies in the interior of Hill's region $\mathfrak{K}_c$. Then $\dot q(t_0)\neq 0$, so $q$ is an immersion near $t_0$. 

{\em Case 2: }$q_0\neq 0$ lies on the boundary of Hill's region $\mathfrak{K}_c$, i.e., the velocity $\dot q(t_0)$ vanishes. 

\begin{figure}[!htb]
\centering
\def\svgwidth{1.0\textwidth}%
\begingroup\endlinechar=-1
\resizebox{1.0\textwidth}{!}{%
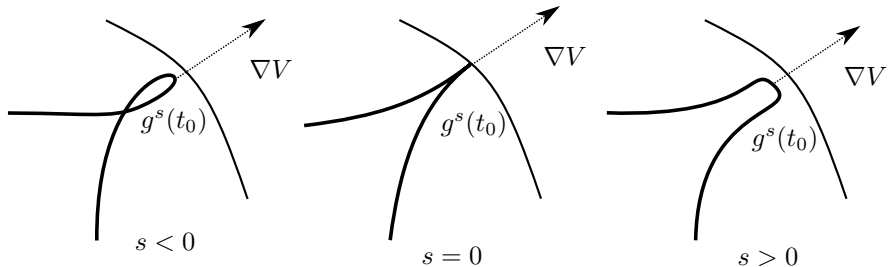%
}\endgroup
\caption{Passing through a cusp at the zero velocity curve (depicted in black)}
\label{fig:cusp-zerovelocity}
\end{figure}

\begin{lemma}\label{lem:cusp-zerovelocity}
Suppose that $\dot q(t_0)=0$ and $B(q_0)\neq 0$. Then the orbit $q$ has a cubical cusp at $t=t_0$. Moreover, in a $1$-parameter family $q^s$ of nearby orbits of the same energy with $q^0=q$ in which the component of $\dot q^s(t_0)$ normal to $\nabla V(q^s(t_0))$ changes sign, an exterior loop shrinks to a cusp and disappears (or the reverse) as shown in Figure~\ref{fig:cusp-zerovelocity}.
%\marginpar{Add figure.}
\end{lemma}

\begin{figure}[!htb]
\centering
\def\svgwidth{0.450\textwidth}%
\begingroup\endlinechar=-1
\resizebox{0.450\textwidth}{!}{%
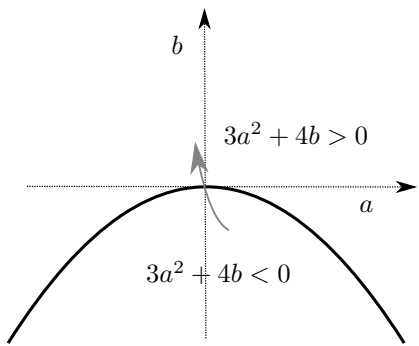%
}\endgroup
\caption{Crossing the discriminant curve}
\label{fig:discriminant}
\end{figure}

{\bf Proof: }
Let us take $t_0=0$ and consider an orbit $q=q^s$ in a $1$-parameter family near $q^0$. Denote by $\dot q_0,\ddot q_0,\dddot q_0$ its derivatives at time $t=0$. 
From Newton's equation $\ddot q=iB(q)\dot q-\nabla V(q)$ and its time derivative we obtain with $B_0=B(q_0)$:
\begin{align*}
   \ddot q_0 &= iB_0\dot q_0-\nabla V(q_0),\cr 
   \dddot q_0 &= iB_0\ddot q_0-D^2V(q_0)\dot q_0 \cr 
   &= -iB_0\nabla V(q_0) - (B_0^2+D^2V(q_0))\dot q_0.
\end{align*}
From the non-vanishing of $\nabla V(q_0)$ and $B_0$ we conclude that for small $\dot q_0$ the vectors $\ddot q_0$ and $\dddot q_0$ are nonzero and roughly orthogonal. Picking complex coordinates in which $q_0=0$, $\ddot q_0=2$ and $\dddot q_0=6i$ and treating the initial velocity $\dot q_0=a+ib$ as a complex parameter, the Taylor expansion of $q$ to 3rd order becomes 
\begin{align*}
   q(t) 
   &= q_0 + \dot q_0t + \frac12\ddot q_0t^2 + \frac16\dddot q_0t^3 + O(t^4) \cr
   &= (at+t^2) + i(bt+t^3) + O(t^4).
\end{align*}
Neglecting the terms of order $4$ or higher, a short computation shows that for parameters satisfying $3a^2+4b<0$ the curve $q$ has a loop which vanishes as the parameters cross the discriminant curve $\{3a^2+4b=0\}$, see Figure~\ref{fig:discriminant}.
%\marginpar{Add figure.} 
Note that a $1$-parameter family in which the component $b$ of $\dot q^s(t_0)$ normal to $\nabla V(q^s(t_0))$ changes sign corresponds to a crossing of the discriminant curve at the origin. 
Since this picture persists under perturbations of order $4$ or higher, the lemma is proved.
\hfill$\square$

{\em Case 3: }$q_0=0$ is a collision. 

\begin{figure}[!htb]
\centering
\def\svgwidth{1.0\textwidth}%
\begingroup\endlinechar=-1
\resizebox{1.0\textwidth}{!}{%
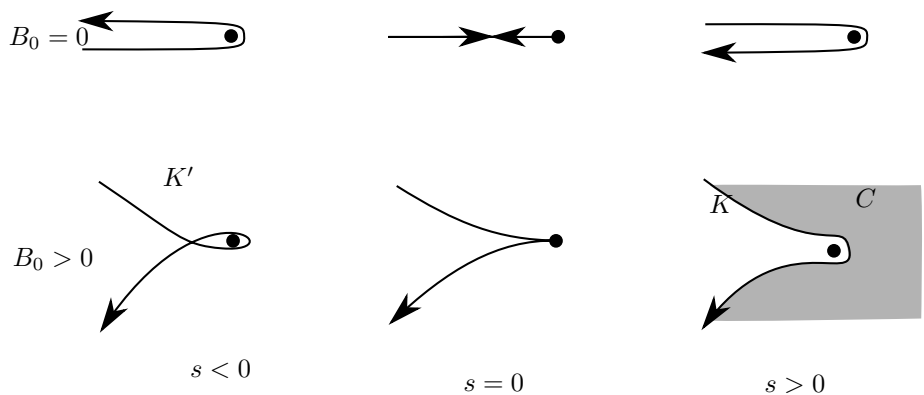%
}\endgroup
\caption{Passing through a cusp at the origin (depicted in black)}
\label{fig:cusp-collision}
\end{figure}

\begin{lemma}\label{lem:cusp-collision}
Suppose that $q_0=0$ and $B(q_0)\neq 0$. Then the orbit $q$ has a cusp at $t=t_0$. Moreover, in a $1$-parameter family $q^s$ of nearby orbits of the same energy with $q^0=q$, a loop around $0$ shrinks to a cusp and becomes a curve going around $0$ the other way (or the reverse) as shown at the bottom of in Figure~\ref{fig:cusp-collision}.
%\marginpar{Add figure.} 
\end{lemma}

{\bf Proof: }
Consider first the situation with $B_0=B(q_0)=0$. Then a collision orbit $q$ bounces back from the origin, and a $1$-parameter family $q^s$ with $q^0=q$ changes the direction in which it encircles the origin as shown at the top of Figure~\ref{fig:cusp-collision}.
If we now switch on a magnetic field with $B_0>0$, then the collision orbit curves everywhere to the left and thus becomes a cusp at the origin. Moreover, the $1$-parameter family $q^s$ now looks as at the bottom of in Figure~\ref{fig:cusp-collision}.
(This holds because away from the origin $q^s$ is close to the cusp $q^0$ for small $s$, whereas near the origin the Newtonian force dominates the magnetic field and thus $q^s$ looks like in the case $B_0=0$.)
\hfill$\square$

\begin{rem}
Performing a time change $\tau=t^{1/3}$ and considering the Taylor expansion of $q(\tau)$ to fourth order near $\tau=0$ one finds that the cusp of $q$ at the origin is {\em quintical}, i.e., equivalent to $\tau\mapsto(\tau^2,\tau^5)$ in suitable coordinates.
\end{rem}

{\bf Stark-Zeeman homotopies. }
Suppose now that $q:\R/T\Z$ is a simple (i.e., not multiply covered) periodic orbit of energy $c$. Since the position and velocity at any time $t_0$ uniquely determine the orbit, it follows that $q$ cannot have any {\em direct self-tangency}, i.e., a self-tangency at which the two velocities agree. By contrast, {\em inverse self-tangencies} at which the two velocities are opposite to each other can occur. Generically, $q$ will be a {\em generic immersion}, i.e.~an immersion with only transverse double points, avoiding the boundary of Hill's region as well as the origin. In view of this discussion and the previous two lemmas, the following definition captures all events that can happen in a generic $1$-parameter family of simple periodic orbits of (varying) planar Stark-Zeeman systems. 

\begin{figure}[!htb]
\centering
\def\svgwidth{0.75\textwidth}%
\begingroup\endlinechar=-1
\resizebox{0.75\textwidth}{!}{%
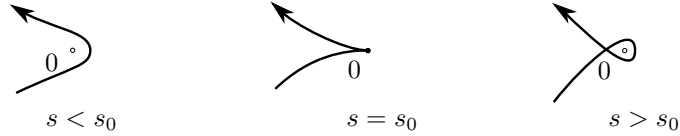%
}\endgroup
\caption{Event ($I_0$): development of interior loop through a cusp}
\label{fig:interior-loop-birth}
\end{figure}

\begin{figure}[!htb]
\centering
\def\svgwidth{0.75\textwidth}%
\begingroup\endlinechar=-1
\resizebox{0.75\textwidth}{!}{%
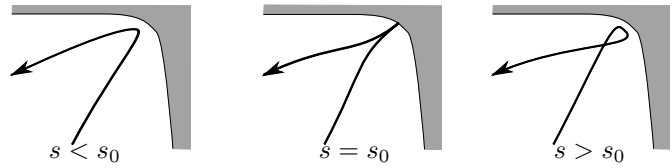%
}\endgroup
\caption{Event ($I_\infty$): development of exterior loop through a cusp}
\label{fig:exterior-loop-birth}
\end{figure}

\begin{figure}[!htb]
\centering
\def\svgwidth{0.75\textwidth}%
\begingroup\endlinechar=-1
\resizebox{0.75\textwidth}{!}{%
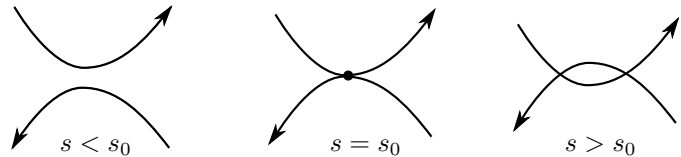%
}\endgroup
\caption{Event ($II^+$), crossing through inverse self-tangency}
\label{fig:inverse-self-tangent}
\end{figure}

\begin{figure}[!htb]
\centering
\def\svgwidth{0.75\textwidth}%
\begingroup\endlinechar=-1
\resizebox{0.75\textwidth}{!}{%
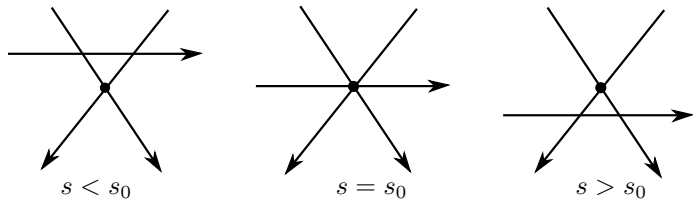%
}\endgroup
\caption{Event ($III$),
an ordinary triple point in a $1$-parameter family}
\label{fig:triple-point}
\end{figure}

\begin{fed}
\label{def:Stark-Zeeman-homotopy}
%\marginpar{Add figure of these four moves.}
A {\em Stark-Zeeman homotopy} is a $1$-parameter family $(K^s)_{s\in[0,1]}$ of closed curves in $\C$ which are generic immersions in $\C^*=\C\setminus\{0\}$ except for the following events at finitely many $s\in(0,1)$: 
\begin{description}
\item[($I_0$)] birth or death of interior loops through cusps at the origin;
\item[($I_\infty$)] birth or death of exterior loops through cusps at infinity;
\item[($II^+$)] crossings through inverse self-tangencies;
\item[($III$)] crossings through triple points.
\end{description}
\end{fed}
These special events, where a $1$-parameter family $K^s$ fails to a generic immersion are illustrated in Figures~\ref{fig:interior-loop-birth}, \ref{fig:exterior-loop-birth}, \ref{fig:inverse-self-tangent} and \ref{fig:triple-point}

Conversely, the following proposition shows that any Stark-Zeeman homotopy can be realized as a family of periodic orbits of Stark-Zeeman systems. 

\begin{prop}\label{prop:SZ-homotopy}
A $1$-parameter family $(K^s)_{s\in[0,1]}$ of closed curves in $\C^*$ that are not multiply covered is a Stark-Zeeman homotopy if and only if there exists a family of diffeomorphisms $F^s:\C^*\to\C^*$ such that the (suitably parametrized) curves $F^s(K^s)$ are simple periodic orbits (possibly with collisions) in a generic family of Stark-Zeeman systems.  
%\marginpar{Be more precise about ``genericity''?}
\end{prop}

\begin{rem}
The restricted three-body problem comes with an anti-symplectic involution $(q_1,q_2;p_1,p_2)\mapsto(q_1,-q_2,-p_1,p_2)$ making this problem very non-generic. For example, if a cusp develops in the upper half-plane, then this cusp has a paired cusp in the lower half-plane. When this occurs with an event of type $I_0$ one obtains a double cusp through $0$.
\end{rem}

The proof of this proposition uses the following lemma. The domain of each curve in the lemma is an open interval.

\begin{lemma}\label{lem:3-curves}
(a) Let $\gamma_1,\gamma_2,\gamma_3$ be three embedded curves through the origin $0$ in $\R^2$ whose derivatives at $0$ are pairwise linearly independent. Let $\wt\gamma_1,\wt\gamma_2, \wt\gamma_3$ be another such triple such that $\wt\gamma_i$ is tangent to $\gamma_i$ at $0$ for $i=1,2,3$. Then there exist neighbourhoods $U$, $U'$ and a diffeomorphism $F:U\to U'$ with $DF(0)=\id$ such that $F(\gamma_i)$ is a reparametrization of $\wt\gamma_i$ for $i=1,2,3$. 

(b) Let $\gamma_1,\gamma_2,\wt\gamma_1,\wt\gamma_2$ be four embedded curves through the origin $0$ in $\R^2$ all having the same nonzero derivative at $0$. Suppose that the corresponding curvatures at $0$ satisfy $\kappa_2-\kappa_1=\wt\kappa_2-\wt\kappa_1\neq 0$. Then there exist neighbourhoods $U$, $U'$ of $0$ and a diffeomorphism $F:U\to U'$ with $DF(0)=\id$ such that $F(\gamma_i)$ is a reparametrization of $\wt\gamma_i$ for $i=1,2$. 
\end{lemma}

{\bf Proof: }
(a) By transitivity, it suffices to prove this for the case that the $\gamma_i$ are straight lines. Moreover, it suffices to consider the case $\gamma_i=\wt\gamma_i$ for $i=1,2$ (because we can then arrange this condition by applying the result twice with some arbitrary auxiliary lines). 
We pick linear coordinates in which $\gamma_1$ and $\wt\gamma_1$ both map to the $x$-axis, $\gamma_2$ and $\wt\gamma_2$ both to the $y$-axis, and $\gamma_3$ to the line $\{y=ax\}$ for some $a\neq 0$. 
Then for a sufficiently small neighbourhood $U'$ of the origin, the image of the curve $\wt\gamma_3$ is the graph $\{y=g(x)\}$ of a smooth function $g$ with $g(0)=0$ and $g'(0)=a\neq 0$. 
By Taylor's formula we can write $g(x)=axh(x)$ for a unique smooth function $h$ with $h(0)=1$. 
By choosing a sufficiently small neighbourhood $U$ as domain we obtain a diffeomorphism $F(x,y):=(x,yh(x))=(X,Y)$ that satisfies $DF(0)=\id$ and preserves the $x$- and $y$-axes. 
Moreover, it sends the graph $y=ax$ to the graph $Y=yh(x)=axh(x)=g(x)=g(X)$ and thus the image of $\gamma_i$ to the image of $\wt\gamma_i$ for $i=1,2,3$.

(b) Arguing as in (a), we see that it suffices to prove this for the case that $\gamma_1=\wt\gamma_1$ is the $x$-axis, $\gamma_2$ is the graph of the function $y=\frac12 ax^2$, and $\wt\gamma_2$ is the graph of a function $y=g(x)$ with $g(0)=g'(0)=0$ and $g''(0)=a\neq 0$. By Taylor's formula we can write $g(x)=\frac12 ax^2h(x)$ for a unique smooth function $h$ with $h(0)=1$. Choosing neighbourhoods as before, we find a diffeomorphism $F(x,y):=(x,yh(x))=(X,Y)$ that satisfies $DF(0)=\id$ and preserves the $x$- and $y$-axes. Moreover, it sends the graph $y=\frac12 ax^2$ to the graph $Y=yh(x)=\frac12 ax^2h(x)=g(x)=g(X)$ and thus the image of $\gamma_i$ to the image of $\wt\gamma_i$ for $i=1,2$.
\hfill$\square$
\medskip

{\bf Proof of Proposition~\ref{prop:SZ-homotopy}: }
As explained above, simple periodic orbits of Stark-Zeeman systems are immersions without direct self-tangencies and with cusps at the zero velocity curves and at the origin. Hence Lemma~\ref{lem:cusp-zerovelocity} and Lemma~\ref{lem:cusp-collision} imply that a generic $1$-parameter family of such orbits is a Stark-Zeeman homotopy, and the latter condition is preserved by composition with diffeomorphisms of the plane. 

Conversely, suppose that $(K^s)_{s\in[0,1]}$ is a Stark-Zeeman homotopy. We will construct a smooth family of Stark-Zeeman systems such that the $K^s$ agree (suitably parametrized) with periodic solutions of energy $0$ of the ODE
\begin{equation}\label{eq:SZ-par}
   \ddot q=B^s(q)i\dot q-\nabla V^s(q). 
\end{equation}
We first pick a smooth family of potentials $V^s$ with a Coulomb singularity at the origin and without critical points such that the curves $K^s$ lie in the interior of the Hill's region $\mathfrak{K}^s=\{-\infty<V^s\leq 0\}$, except for their cusps which lie on the boundary $\{0\}\cup\{V^s=0\}$. We add preliminary magnetic fields $B^s_0$ near the cusps and apply diffeomorphisms there to make $K^s$ agree with zero energy solutions of~\eqref{eq:SZ-par} near the cusps. We cut off $B^s_0$ and extend this function by zero outside neighbourhoods of the cusps. 
Next, we use the neighbourhoods and diffeomorphisms provided by Lemma~\ref{lem:3-curves} near all double and triple points to map $K^s$ to solutions of~\eqref{eq:SZ-par} near double and triple points. These diffeomorphisms can be chosen smoothly in $s$, agreeing with the given ones near the cusps, and can be extended to a smooth family of diffeomorphisms $F^s$. 
We parametrize the curves $F^s(K^s)$ near the cusps and double/triple points as solutions of~\eqref{eq:SZ-par}, and extend these parametrizations to parametrizations $q^s:[0,T^s]\to\C^*$ of $F^s(K^s)$ by the energy requirement $\frac12|\dot q^s|^2+V^s(q^s)=0$. 
Note that $\dot q^s\neq 0$ away from the cusps.
Taking a time derivative of the energy condition we find $0=\la\dot q^s,\ddot q^s+\nabla V^s(q^s)\ra$, so $\ddot q^s+\nabla V^s(q^s)$ is a multiple of $i\dot q^s$ away from the cusps. 
Since $q^s$ is injective away from the double and triple points, we thus find unique well-defined functions $B^s(q)$ along $F^s(K^s)$ such that $\ddot q^s+\nabla V^s(q^s) = i B^s(q^s)\dot q^s$, i.e., $q^s$ is a zero energy solution of~\eqref{eq:SZ-par}. 
We arbitrarily extend $B^s$ away from $K^s$ in a way that smoothly depends on $s$ and Proposition~\ref{prop:SZ-homotopy} is proved.
\hfill$\square$

%%%
\subsection{Topology of plane curves}
%%%

In this subsection we recall some facts from~\cite{arnold} about immersed loops in the plane, which will be applied to Stark-Zeeman systems in the next subsection. In particular, we recall the definition of Arnold's $J^+$-invariant and derive some of its properties. 

\begin{figure}[!htb]
\centering
\def\svgwidth{0.75\textwidth}%
\begingroup\endlinechar=-1
\resizebox{0.75\textwidth}{!}{%
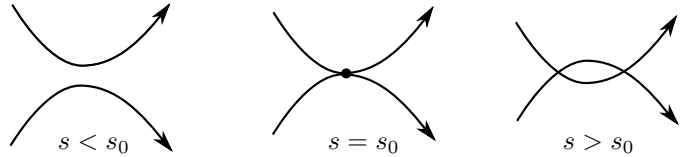%
}\endgroup
\caption{A $1$-parameter family of immersions going through a direct self-tangency}
\label{fig:direct-self-tangent}
\end{figure}

\begin{figure}[!htb]
\centering
\def\svgwidth{1.0\textwidth}%
\begingroup\endlinechar=-1
\resizebox{1.0\textwidth}{!}{%
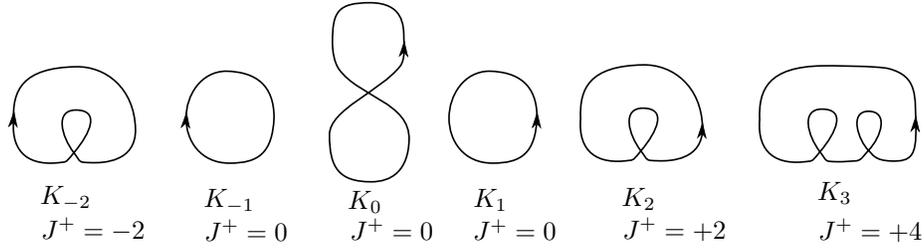%
}\endgroup
\caption{The standard curves and their $J^+$.}
\label{fig:standard-curves}
\end{figure}

By an {\em immersion} we will mean an immersion $q:S^1\to\C$ of the circle into the complex plane, considered up to orientation preserving reparametrization. Note that $q$ is canonically oriented (by $\dot q$) and cooriented (by $i\dot q$). By a slight abuse of notation we will identify an immersion with its image $K\subset\C$. A generic immersion has only transverse double points. During a generic homotopy of immersions the following three events occur: {\em triple points} (as in Figure~\ref{fig:triple-point}), {\em direct self-tangencies} (where the coorientations agree, see Figure~\ref{fig:direct-self-tangent}), and {\em inverse self-tangencies} (where the coorientations do not agree, see Figure~\ref{fig:inverse-self-tangent}). 
The discriminant hypersurfaces corresponding to these events are naturally cooriented, and Arnold defines three invariants (St, $J^+$, $J^-$) of generic immersions in terms of their changes under positive crossings of the discriminant hypersurfaces. They are independent of the orientation and coorientation and additive under connected sums. Of these, the invariant $J^+$ is relevant for us. It increases by $2$ under positive crossings through direct self-tangencies (where the number of double points increases) and remains unchanged under crossings through triple points and inverse self-tangencies. It is normalized to $0$ on the figure eight $K_0$ and to $2-2|j|$ on the circle $K_j$ with $|j|-1$ interior loops and rotation number (winding number of the tangent vector) $j\neq 0$; see Figure~\ref{fig:standard-curves} for pictures of the {\em standard curves} $K_j$. 
%\marginpar{Add figure.}

\begin{figure}[!htb]
\centering
\def\svgwidth{0.40\textwidth}%
\begingroup\endlinechar=-1
\resizebox{0.40\textwidth}{!}{%
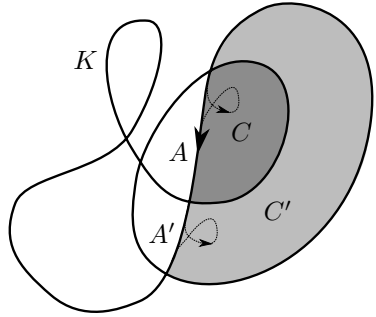%
}\endgroup
\caption{Adding a loop to an arc $A$ of $K$; here $w(K,C)=2$ and $w(K,C')=1$.}
\label{fig:loop}
\end{figure}

\begin{lemma}\label{loop}
For an immersion $K\subset\C$, consider a connected component $C$ of $\C\setminus K$ and a boundary arc $A$ of $C$. Orient $K$ such that $A$ inherits the boundary orientation from $C$ and denote by $w(K,C)\in\Z$ the winding number of $K$ (with this orientation) around a point in $C$. Then the invariant $J^+$ changes by $-2w(K,C)$ under addition of a small loop in $C$ to the arc $A$, see Figure~\ref{fig:loop}. 
%\marginpar{Add figure.}
\end{lemma}

\begin{figure}[!htb]
\centering
\def\svgwidth{0.80\textwidth}%
\begingroup\endlinechar=-1
\resizebox{0.80\textwidth}{!}{%
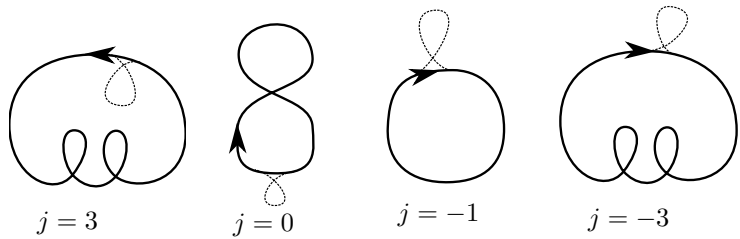%
}\endgroup
\caption{Adding a loop to the standard curves.}
\label{fig:loop-standard-curves}
\end{figure}

\begin{figure}[!htb]
\centering
\def\svgwidth{0.90\textwidth}%
\begingroup\endlinechar=-1
\resizebox{0.90\textwidth}{!}{%
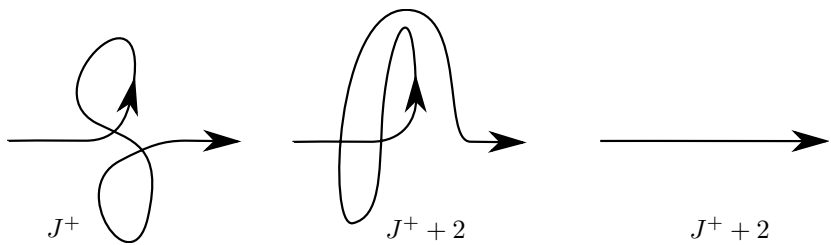%
}\endgroup
\caption{Cancelling two loops changes $J^+$ by $2$.}
\label{fig:cancelling-loops}
\end{figure}

{\bf Proof: }
For $K,C,A$ as in the lemma, let $\wt K$ be the immersion obtained by adding a small loop in $C$ to the arc $A$ of $K$ and set
$$
   I(K,C,A) := J^+(\wt K)-J^+(K)+2w(K,C).
$$ 
Consider a component $C'$ of $\C\setminus K$ adjacent to $C$ with boundary arc $A'$ adjacent to $A$ and winding number $w(K,C')=w(K,C)\pm 1$ as in Figure~\ref{fig:loop}. Let $\wt K'$ be the immersion obtained by adding a small loop in $C'$ to the arc $A'$ of $K$. Figure~\ref{fig:loop} shows that $J^+(\wt K')-J^+(\wt K)=\mp 2$, and therefore $I(K,C,A)=I(K,C',A')$. Thus $I(K,C,A)$ is independent of $C$ and $A$ and we can simply denote it by $I(K)$ (note that at this point $I(K)$ may still depend on the orientation of $K$). It follows that $I(K)$ is invariant under crossings through triple points or (direct or inverse) self-tangencies (because we can always attach the new loop along an arc that is not involved in the crossing, so that $J^+(\wt K)-J^+(K)$ and $w(K,C)$ both do not change under the crossing). Hence the invariant $I$ is uniquely determined by its value on the standard curves $K_j$, which is computed using Figure~\ref{fig:loop-standard-curves} as follows.
%\marginpar{Add figure.}
For $j>0$ we take the component $C$ with $w(K_j,C)=1$; then $\wt K_j=K_{j+1}$ and thus
$$
   I(K_j)=J^+(K_{j+1})-J^+(K_j)+2w(K_j,C) = -2j - (2-2j) + 2 = 0.  
$$
For $j=0$ we take the unbounded component $C$ so that $w(K_0,C)=0$; then $\wt K_0$ is obtained from $K_{1}$ without crossing through a direct self-tangency and thus
$$
   I(K_0)=J^+(K_{1})-J^+(K_0)+2w(K_0,C) = 0 - 0 + 0 = 0.  
$$
For $j=-1$ we take the unbounded component $C$ so that $w(K_{-1},C)=0$; then $\wt K_{-1}=K_{0}$ and thus
$$
   I(K_{-1})=J^+(K_{0})-J^+(K_{-1})+2w(K_0,C) = 0 - 0 + 0 = 0.  
$$
For $j<-1$ we take the unbounded component $C$ so that $w(K_j,C)=0$; Figure~\ref{fig:cancelling-loops} 
%\marginpar{Add figure.} 
shows that $\wt K_j$ is obtained from $K_{j+1}$ by crossing once negatively through a direct double point and thus
$$
   I(K_j)=J^+(K_{j+1})-2-J^+(K_j)+2w(K_j,C) = (4+2j) - 2 - (2+2j) + 0 = 0.  
$$
This shows that $I(K)$ vanishes identically, which proves the lemma. 
\hfill$\square$

In particular, we have the following corollary (which also follows directly from the additivity of $J^+$ under connected sum and the fact that adding an exterior loop corresponds to connected sum with $K_0$). 

\begin{cor}\label{loopcor}
The invariant $J^+$ does not change under addition of an {\em exterior loop}, i.e., a loop in the unbounded component. 
\hfill$\square$
\end{cor}

\subsection{The two invariants $\J_1$ and $\J_2$}
%%%
For a generic immersion $K$ in $\mathbb{C}^*$ denote by $J^+(K) \in 2 \mathbb{Z}$ its
$J^+$-invariant and by $w_0(K) \in \mathbb{Z}$ its winding number around the origin.  We set
$$\mathcal{J}_1(K):=J^+(K)+\frac{w_0(K)^2}{2}.$$
\begin{prop}\label{prop:J1}
$\mathcal{J}_1$ is invariant under Stark-Zeeman homotopies. 
\end{prop}
\textbf{Proof: }We need to show that $\mathcal{J}_1$ is invariant under the moves in Definition~\ref{def:Stark-Zeeman-homotopy}. Under the moves $(II^+)$ and $(III)$ both quantities $J^+$ and $w_0$ are invariant. The same
is true for the move $(I_\infty)$ due to Corollary~\ref{loopcor}. Therefore, $\mathcal{J}_1$ is invariant
under these three moves. To discuss the invariance of $\mathcal{J}_1$ under the move $(I_0)$ first
observe that $J^+$ is invariant under orientation reversion while $w_0$ changes sign, so $\mathcal{J}_1$ is invariant under orientation reversion. Hence we can choose the orientations of the immersions $K$ and $K'$ before and after the move $(I_0)$ as shown at the bottom of Figure~\ref{fig:cusp-collision}. Then their winding numbers are related by $w_0(K')=w_0(K)+2$.
Moreover, $K'$ is obtained from $K$ by adding a small loop in the component $C$ of $\C\setminus K$ shown in the figure (having the arc of $K$ as a positive boundary arc). Since the winding number of $K$ around $C$ is $w_0(K)+1$, 
by Lemma~\ref{loop} the $J^+$-invariants are related by
$$J^+(K')=J^+(K)-2w(K,C)=J^+(K)-2w_0(K)-2.$$
Hence we compute
\begin{eqnarray*}
\mathcal{J}_1(K')&=&J^+(K')+\frac{w_0(K')^2}{2}\\
&=&J^+(K)-2w_0(K)-2+\frac{(w_0(K)+2)^2}{2}\\
&=&J^+(K)-2w_0(K)-2+\frac{w_0(K)^2+4w_0(K)+4}{2}\\
&=&J^+(K)+\frac{w_0(k)^2}{2}\\
&=&\mathcal{J}_1(K).
\end{eqnarray*}
This proves invariance of $\mathcal{J}_1$ under the move $(I_0)$ and hence the proposition. \hfill $\square$
\\ \\
Recall from Section~\ref{levicivita} the Levi-Civita regularization map
$$L \colon \mathbb{C}^* \to \mathbb{C}^*, \quad v \mapsto v^2$$
which regularizes collisions at the origin. Note that the preimage $L^{-1}(K)$ of a generic immersion $K\subset\C^*$ is again a generic immersion in $\C^*$. We distinguish two cases.

{\em Case 1: }$w_0(K)$ is odd. Then the preimage $L^{-1}(K)$ is connected and we set
$$\mathcal{J}_2(K):=J^+(L^{-1}(K)).$$
Note that $L:L^{-1}(K)\to K$ is a $2$--$1$ covering, in particular $w_0(L^{-1}(K))=2w_0(K)$.  

{\em Case 2: }$w_0(K)$ is even. Then the preimage $L^{-1}(K)$ is disconnected and we set
$$\mathcal{J}_2(K):=J^+(\wt K)$$ 
for one component $\wt K$ of $L^{-1}(K)$. Since the two components of $L^{-1}(K)$ are related by a $180$ degree rotation, this definition does not depend on the choice of $\wt K$. 
Note that $L:\wt K\to K$ is a diffeomorphism, in particular $w_0(\wt K)=w_0(K)$.  

\begin{prop}\label{prop:J2}
$\mathcal{J}_2$ is invariant under Stark-Zeeman homotopies.
\end{prop}
\textbf{Proof: }Since the Levi-Civita regularization map regularizes collisions, the quantity $\mathcal{J}_2$
is invariant under the move $(I_0)$. By definition of the $J^+$-invariant, it is invariant as well under the moves
$(II^+)$ and $(III)$ and invariance under the move $(I_\infty)$ is the content of Corollary~\ref{loopcor}.
This proves the proposition. \hfill $\square$

\medskip
Propositions~\ref{prop:J1} and~\ref{prop:J2} together prove Theorem A of the Introduction.

%%%
\subsection{Comparison of the two invariants}
%%%
Next we will discuss the relation between the invariants $\mathcal{J}_1$ and $\mathcal{J}_2$. We begin with the case that $w_0(K)$ is odd.

\begin{prop}\label{prop:comparison-odd}
If $w_0(K)$ is odd, then $\mathcal{J}_2(K)=2\mathcal{J}_1(K)-1$. 
\end{prop}

\textbf{Proof: }
On generic immersions $K\subset\C^*$ with $w_0(K)$ odd consider the quantity $I(K):=\J_2(K)-2\J_1(K)$.
Since $L:\C^*\to\C^*$ is an oriented $2$--$1$ covering, each point $p$ of $K$ has two preimages in $L^{-1}(K)$ near which $L^{-1}(K)$ looks like $K$ near $p$. This shows that $I$ is invariant under the moves $(II^\pm)$ and $(III)$. Using these moves and adding exterior loops to adjust the rotation number (which by Corollary~\ref{loopcor} does not affect $I$), we can deform $K$ to any given generic immersion with the same winding number around the origin without changing $I$. So it remains to compute $I(K^w)$ on one generic immersion $K^w$ of winding number $w_0(K^w)=w$ for each odd integer $w$. By invariance of $I$ under orientation reversal, it suffices to consider positive $w$. 

\begin{figure}[!htb]
\centering
\def\svgwidth{0.70\textwidth}%
\begingroup\endlinechar=-1
\resizebox{0.70\textwidth}{!}{%
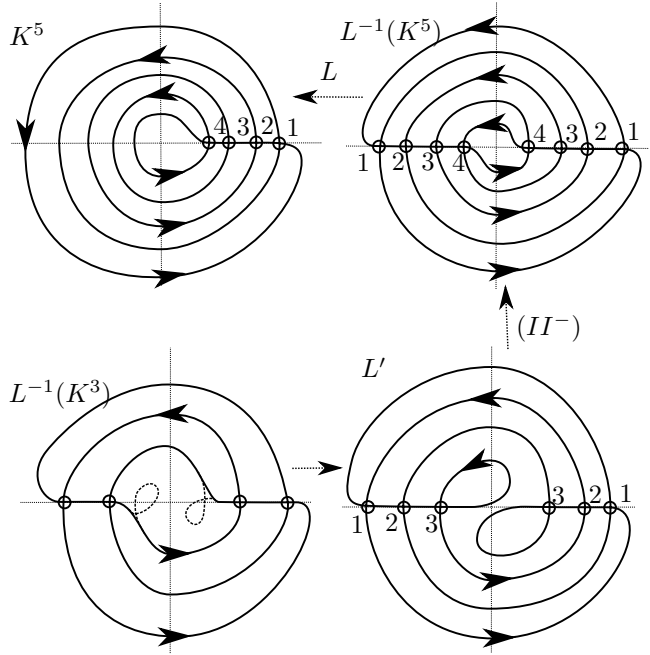%
}\endgroup
\caption{The curve $K^5$ and its Levi-Civita lift $L^{-1}(K^5)$.}
\label{fig:J2}
\end{figure}

For any $w\in\N$ consider the immersion $K^w$ of winding number $w_0(K^w)=w$ as shown at the top left in Figure~\ref{fig:J2}.
%\marginpar{Add figure}

We see that $K^w$ is obtained from $K^{w-1}$ by adding a small loop in its innermost component $C$, which has winding number $w-1$ around the origin. In view of Lemma~\ref{loop}, this gives the recursion relation $J^+(K^w)=J^+(K^{w-1})-2(w-1)$. Together with the initial condition $J^+(K^1)=0$ (which holds because $K^1$ equals the standard curve $K_1$), this determines $J^+(K^w)$ to be
\begin{equation*}
   J^+(K^w) = -2(1+2+3+\cdots+(w-1)) = -w(w-1). 
\end{equation*}
Now suppose that $w\in\N$ is odd and consider the immersion $L^{-1}(K^w)$ as depicted at the top right in Figure~\ref{fig:J2}. Applying a $(II^-)$ move (which does not change $J^+$) to the two innermost arcs yields the immersion $L'$ shown at the bottom right in Figure~\ref{fig:J2}.
The figure shows that $L'$ is obtained from $L^{-1}(K^{w-2})$ by adding two disjoint small loops in its innermost component $C$, which has winding number $w-2$ around the origin. In view of Lemma~\ref{loop}, this gives the recursion relation $J^+(L^{-1}(K^w))=J^+(L')=J^+(L^{-1}(K^{w-2}))-4(w-2)$. Together with the initial condition $J^+(L^{-1}(K^1))=0$ (which holds because $L^{-1}(K^1)$ equals the standard curve $K_1$), this determines $J^+(L^{-1}(K^w))$ to be
\begin{equation*}
   J^+(L^{-1}(K^w)) = -4(1+3+5+\cdots+(w-2)) = -(w-1)^2. 
\end{equation*}
Combining the preceding two displayed equations, it follows that
$$
   I(K^w) = J^+(L^{-1}(K^w)) - 2\bigl(J^+(K^w) + w^2/2\bigr) = -(w-1)^2 + 2w(w-1) - w^2 = -1.
$$
In view of the discussion at the beginning of the proof, this shows that for each generic immersion $K\subset\C^*$ with odd winding number around the origin we have $I(K)=\mathcal{J}_2(K)-2\mathcal{J}_1(K)=-1$. 
\hfill $\square$
\\ \\ 

\begin{figure}[!htb]
\centering
\def\svgwidth{0.60\textwidth}%
\begingroup\endlinechar=-1
\resizebox{0.60\textwidth}{!}{%
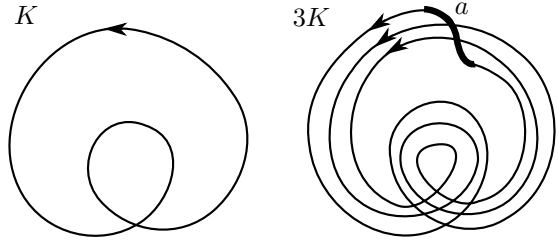%
}\endgroup
\caption{The $3$-satellite $3K$.}
\label{fig:satellite}
\end{figure}

{\bf Satellites. }
Next we turn to the case that $w_0(K)$ is even. We will generate examples using the following construction. For an oriented generic immersion $K\subset\C^*$ and $n\in\N$, the {\em $n$-satellite} $nK$ is defined as follows, see Figure~\ref{fig:satellite}.
%\marginpar{Add figure}
We start at some point $p\in K$ and construct a spiral along $K$ slowly moving off linearly to the left of $K$. After $n$ turns along $K$, we close off the spiral by a short arc $a$ intersecting it $n-1$ times. Thus $nK$ has $n^2$ double points for each double point of $K$, plus the $n-1$ double points on the arc $a$. The rotation number of $nK$ is $n$ times that of $K$, and $w_0(nK)=nw_0(K)$. Changing the orientation of $K$ leads to an $n$-satellite which differs from that of $K$ by a reflection of the plane and thus has the same invariants $\J_1$ and $\J_2$.

Note that $nK_{-1}$ (with $K_{-1}$ enclosing the origin) equals the immersion $K^n$ with $w_0(K^n)=n$ defined above. 

\begin{lemma}\label{2-satellite}
For $j\in\N$, the $2$-satellite $2K_j$ of the standard curve $K_j$ with $w_0(K_j)=1$ has the invariants 
$$
   \J_1(2K_j) = 2\J_2(2K_j) = -8(j-1).
$$
More generally, for any generic immersion $K\subset\C^*$ with $w_0(K)=1$, the invariants of its $2$-satellite $2K$ are related by $\J_1(2K)=2\J_2(2K)$.   
\end{lemma}

\begin{figure}[!htb]
\centering
\def\svgwidth{0.750\textwidth}%
\begingroup\endlinechar=-1
\resizebox{0.750\textwidth}{!}{%
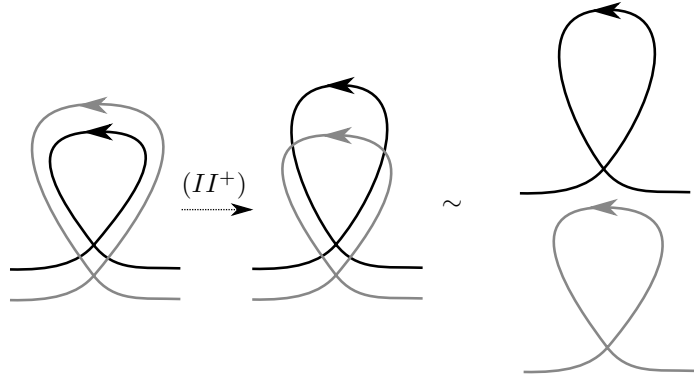%
}\endgroup
\caption{Resolving a double loop.}
\label{fig:double-loop}
\end{figure}

{\bf Proof: }
Consider the double loop in $2K_j$ corresponding to a loop in $K_j$ as shown on the left in Figure~\ref{fig:double-loop}.
%\marginpar{Add figure}
Pulling the interior loop down by a $(II^+)$ move (followed by a $(III)$ and two $(II^-)$ moves) yields an immersion $K'$ with $J^+(K')=J^+(2K_j)+2$ in which the double loop has been turned into two single loops as shown on the right in Figure~\ref{fig:double-loop}.
Removing the resulting two loops from $K'$ yields the satellite immersion $2K_{j-1}$. By Lemma~\ref{loop} the invariants are related by $J^+(K')=J^+(2K_{j-1})-2-4$, so we get the recursion relation $J^+(2K_j)=J^+(2K_{j-1})-8$. Together with the initial condition $J^+(2K_1)=-2$ (since $2K_1=K_2$), this shows $J^+(2K_j)=-2-8(j-1)$ and therefore $\J_1(2K_j)=J^+(2K_j)+2^2/2 = -8(j-1)$. 

For $\J_2$, note that for each double loop in $2K_j$ the intersection points between the two loops disappear in the lift $\wt{2K_j}$ under the Levi-Civita covering. Hence $\wt{2K_j}$ has $2(j-1)$ disjoint loops, so it coincides with the standard immersion $K_{2j-1}$ and we find $\J_2(2K_j)=J^+(K_{2j-1})=-4(j-1)$. 

Finally, consider the quantity $I(K):=\J_1(2K)-2\J_2(2K)$ on generic immersions $K\subset\C^*$ with $w_0(K)=1$. Since every double point in $K$ gives rise to $4$ double points in $2K$ and $2$ double points in its lift $\wt{2K}$ (all of the same type), it follows that $I$ is invariant under the moves $(II^\pm)$ and $(III)$. Using these moves and adding exterior loops (which by Corollary~\ref{loopcor} does not affect $I$), we can deform any given $K$ to some $K_j$ and conclude $I(K)=I(K_j)=0$. 
This proves the lemma. \hfill$\square$
\\ \\

The following proposition shows that for even winding number around the origin the invariants $\J_1$ and $\J_2$ are completely independent.

\begin{prop}\label{prop:comparison-even}
On generic immersions $K\subset\C^*$ with $w_0(K)$ even, the pair of invariants $(\J_1(K),\J_2(K))$ attains all values in $2\Z\times 2\Z$. 
\end{prop}

{\bf Proof: }
We start with the $2$-satellite $2K_j$ of a standard knot $K_j$ with $w_0(K_j)=1$ and $j\in\N$, which by Lemma~\ref{2-satellite} has invariants 
$$\J_1(2K_j) = -8(j-1)\quad\text{and}\quad \J_2(2K_j) = -4(j-1).$$ 
A $(II^+)$ move pulling one strand of $2K_j$ across the neighbouring strand increases $\J_1(2K_j)$ by $2$ and leaves $\J_2(2K_j)$ unchanged (because the two new double points do not give rise to double points in the lift $\wt{2K_j}$ under the Levi-Civita covering). Performing $k\in\N_0$ such operations, we obtain an immersion $K_{j,k}$ with invariants 
$$\J_1(K_{j,k}) = -8(j-1)+2k\quad\text{and}\quad \J_2(K_{j,k}) = -4(j-1).$$ 
Finally, we take the connected sum $K_{j,k,\ell}$ of $K_{j,k}$ and an immersion $L$ with $w_0(L)=0$ and $J^+(L)=2\ell$, for any $\ell\in\Z$. Its lift $\wt{K_{j,k,\ell}}$ under the Levi-Civita covering is the connected sum of $\wt{K_{j,k}}$ and $L$, so by additivity of $L^+$ we get the invariants 
\begin{equation}\label{eq:invariants}
\J_1(K_{j,k,\ell}) = -8(j-1)+2k+2\ell\quad\text{and}\quad \J_2(K_{j,k,\ell}) = -4(j-1)+2\ell.
\end{equation} 
By appropriate choices of $j\in\N$, $k\in\N_0$ and $\ell\in\Z$ we can arrange arbitrary values in $2\Z\times 2\Z$ for the pair $\bigl(\J_1(K_{j,k,\ell}),\J_2(K_{j,k,\ell})\bigr)$. This proves the proposition for winding number $2$ around the origin.

For arbitrary even winding number $2w$, let $K^w$ be any generic immersion with $w_0(K^w)=w$. Its $2$-satellite $2K^w$ has invariants $\bigl(\J_1(2K^w),\J_2(2K^w)\bigr)=(2a,2b)$ for some $a,b,\in\Z$ (depending on $w$). Adding a loop to $K^w$ in a component of winding number $1$ corresponds to adding a double loop to $2K^w$, which becomes two disjoint loops in the lift $\wt{2K^w}$. Hence this operation decreases $\J_2(2K^w)$ by $4$, and thus $\J_1(2K^w)$ by $8$ according to Lemma~\ref{2-satellite}. Adding $(j-1)$ such loops thus decreases $\J_1(2K^w)$ by $8(j-1)$ and $\J_2(2K^w)$ by $4(j-1)$. Performing subsequently the other two operations described above with $k\in\N_0$ and $\ell\in\Z$, we can therefore change the original values $(2a,2b)$ by arbitrary pairs of numbers as in equation~\eqref{eq:invariants}, and thus achieve all values in $2\Z\times 2\Z$. 
\hfill$\square$

\medskip
Propositions~\ref{prop:comparison-odd} and~\ref{prop:comparison-even} together prove Theorem B of the Introduction.

%%%
\subsection{Further discussion}
%%%
{\bf Smooth knot type. }
Every periodic orbit $\gamma$ of energy $c<c_1$ in a planar Stark-Zeeman system describes a smooth knot in the Moser regularized energy hypersurface $\mathcal{S}_c^b\cong\R P^3$. If $\gamma$ avoids the zero velocity locus and collisions, then it is obtained from its footpoint projection $K\subset\R^2$ by adding to each point on $K$ its velocity vector. Clearly, the knot type of $\gamma$ is an invariant of $K$ under Stark-Zeeman homotopies. 

Alternatively, we can associate to each point on $K$ its (appropriately scaled) normal vector; the resulting knot in $\R P^3$ differs from $\gamma$ just by a rotation by $-\pi/2$ and thus represents the same knot type.  
According to~\cite[Theorem 1.2]{chmutov-goryunov-murakami}, every oriented knot type in $\R^2\times S^1$ is obtained by adding the normal vectors to an immersion in the plane. Since every periodic orbit in $\R P^3$ can be pushed into $\R^2\times S^1$, the complement of the zero velocity and collision loci, this implies 
\begin{cor}
Every oriented knot type in $\R P^3$ is realized by a Moser regularized periodic orbit in some planar Stark-Zeeman system.\hfill$\square$
\end{cor}
We can further lift the orbit $\gamma$ from $\R P^3$ to an orbit $\wt\gamma$ in the Levi-Civita regularization $S^3$ (choosing one of the two components if $\gamma$ is contractible and traversing $\gamma$ twice if it is noncontractible). The oriented knot type of $\wt\gamma$ is then also an invariant of the footpoint projection $K\subset\R^2$ under Stark-Zeeman homotopies.  

{\bf Relation to Legendrian knots. }
Adding the unit conormal vectors to an immersion $K\subset\R^2$ actually yields a {\em Legendrian} knot in the unit cotangent bundle $S^*\R^2=\R^2\times S^1$ with its canonical contact structure. The invariant $J^+(K)$ can be defined as the Bennequin-Tabachnikov invariant of its Legendrian lift~\cite{arnold}, and more refined invariants of immersions in the plane can be obtained from Chekanov type invariants of their Legendrian lifts~\cite{ng}. Since the moves ($I_0$) and ($I_\infty$) in Definition~\ref{def:Stark-Zeeman-homotopy} do not lift to Legendrian isotopies, general invariants of Legendrian lifts will not be invariant under Stark-Zeeman homotopies. It would be interesting, however, to explore the possibility of defining more refined invariants by suitably modifying Legendrian isotopy invariants.

{\bf Higher energy levels. }
In the planar restricted three-body system, only the bounded components of energy levels below the first critical value $c_1$ fit our definition of a planar Stark-Zeeman system. For the bounded component of a higher energy level, the potential will have two Coulomb singularities in the Hill's region. The extension of the invariant $\J_1$ to this situation is straightforward, adding correction terms for both winding numbers around the singularities (and similarly for more that two singularities).   For the extension of $\J_2$ we have infinitely many possibilities due to the different Levi-Civita type regularizations corresponding to the different covers of the Moser regularized energy hypersurface $\R P^3\#\R P^3$.  

{\bf Bifurcations. }
The notion of a Stark-Zeeman homotopy captures the events occurring in a generic $1$-parameter family of simple periodic orbits in planar Stark-Zeeman systems. Now such a family generically also undergoes bifurcations, in which an $n$-fold multiple of the orbits becomes degenerate and gives birth to a new family of simple orbits~\cite{abraham-marsden}. The footpoint projection $L$ of one of these new simple orbits is thus a perturbation of the $n$-fold multiple of the projection $K$ of one of the original orbits. Based on the above discussion of satellites, we conjecture that $J^+(L)$ is bounded from below by a constant depending only on $J^+(K)$ and $n$. Such an estimate would be useful for comparing the values of the invariants $\J_1$ and $\J_2$ in different planar Stark-Zeeman systems. 
\\ \\

\end{document}